\documentclass[final]{siamltex}

\usepackage{amssymb} \usepackage{cite} \usepackage{graphicx}
\usepackage{psfrag} \usepackage{subfigure} \usepackage{url}
\usepackage{amsmath}
\usepackage{color}

\usepackage[colorlinks=false,urlcolor=blue,citecolor=blue,linkcolor=blue,
  bookmarks=true, bookmarksopen=false,pdftitle=created by
  dvipdf,pdfcreator=NM,
  pdfauthor=Megiddo,pdfsubject=mor.sty:6.29.04]{hyperref}

\newcommand{\Prob}[1]{\mathsf{Prob\,}#1} \newcommand{\remove}[1]{}

\newcommand{\EXP}[1]{\mathsf{E}\!\left[#1\right] }

\newcommand{\ram}[1]{#1}
\newcommand{\rem}[1]{}
\newcommand{\di}[1]{\mathsf{dist}\left(#1,X^*\right)} \def\rn{\Re^n}
\def\e{\epsilon} \def\l{\left} \def\r{\right} \def\a{\alpha}
\def\m{\mu}  \def\s{\sigma} 
 \def\Argmin{\mathop{\rm Argmin}}
 \def\Argmax{\mathop{\rm Argmax}}
\newtheorem{assumption}{Assumption}

\title{Incremental Stochastic Subgradient Algorithms for Convex
Optimization}

\author{S.~Sundhar~Ram, A.~Nedi\'{c}, and V.~V.~Veeravalli \thanks{The
        first and the third authors are with the Dept.\ of Electrical
        and Computer Eng., University of Illinois at Urbana-Champaign.
        The second author is with the Dept.\ of Industrial and
        Enterprise Systems Eng., University of Illinois at
        Urbana-Champaign. They can be contacted at {\texttt
        \{ssriniv5,angelia,vvv\}@illinois.edu}. }  This work is
        partially supported by Vodafone and the National Science
        Foundation under CAREER grant CMMI 07-42538.  }

\begin{document}
\maketitle

\begin{abstract}
In this paper we study the effect of stochastic errors on two
constrained incremental sub-gradient algorithms. We view the
incremental sub-gradient algorithms as decentralized network
optimization algorithms as applied to minimize a sum of functions,
when each component function is known only to a particular agent of a
distributed network. We first study the standard cyclic incremental
sub-gradient algorithm in which the agents form a ring structure and
pass the iterate in a cycle. We consider the method with stochastic
errors in the sub-gradient evaluations and provide sufficient
conditions on the moments of the stochastic errors that guarantee
almost sure convergence when a diminishing step-size is used. We also
obtain almost sure bounds on the algorithm's performance when a
constant step-size is used.  We then consider \ram{the} Markov
randomized incremental subgradient method, which is a non-cyclic
version of the incremental algorithm where the sequence of computing
agents is modeled as a time non-homogeneous Markov chain. Such a model
is appropriate for mobile networks, as the network topology changes
across time in these networks.  We establish the convergence results
and error bounds for the Markov randomized method in the presence of
stochastic errors for diminishing and constant step-sizes,
respectively.
\end{abstract}
 
\pagestyle{myheadings} \thispagestyle{plain} \markboth{\ \ \
  S.~Sundhar~Ram, A.~Nedi\'{c} and V.~V.~Veeravalli}{Incremental
  Stochastic Subgradient Algorithms \ \ \ }

\section{Introduction}
A problem of recent interest in distributed networks is the design of
decentralized algorithms to minimize a sum of functions, when each
component function is known only to a particular agent
\cite{Nedic01,Rabbat05,Blatt07,Nedic08a, Johansson07}.  Such problems
arise in many network applications, including in-network estimation,
learning, signal processing and resource allocation
\cite{Eryilmaz06,Low99,Predd06,Rabbat05,Sundhar07b,Sundhar07c}.  In
these applications, there is no central coordinator that has access to
all the information and, thus, decentralized algorithms are needed to
solve the problems.  In this paper, we consider decentralized
subgradient methods for constrained minimization of a sum of convex
functions, when each component function is only known partially (with
stochastic errors) to a specific network agent.  We study two
incremental subgradient methods with stochastic errors: a cyclic and a
(non-cyclic) Markov randomized incremental method.

The cyclic incremental algorithm is a decentralized method in which
the network agents form a ring and process the information cyclically.
The incremental method was originally proposed by Kibardin
\cite{Kibardin80} and has been extensively studied more recently in
\cite{Solodov98,Luo94,Gaivoronski94,Nedic01,Blatt07}.  Incremental
gradient algorithms were first used for optimizing the weights in
neural network training \cite{Luo94,Solodov98,Gaivoronski94}, and most
of the associated literature deals with differentiable non-convex
unconstrained optimization problems
\cite{Solodov98,Solodov98b,Gaivoronski94,Bertsekas00,Blatt07}.  The
incremental subgradient algorithm for non-differentiable constrained
convex optimization has been investigated in \cite{Nedic01,Nedic01b}
without errors, and in \cite{Kiwiel03,Solodov98b,Nedic07, Rabbat05}
where the effects of deterministic errors are considered.  The
algorithm that we consider in this paper is stochastic and as such
differs from the existing literature.  For this algorithm, we
establish convergence for diminishing step-size and provide an error
bound for constant step-size.

The Markov randomized incremental algorithm is a decentralized method
where the iterates are generated incrementally within the network by
passing them from agent to agent.  Unlike the cyclic incremental
method, where the agent network has a ring structure and the
information flow is along this ring (cycle), in the Markov randomized
incremental method, the network can have arbitrary structure.
However, similar to cyclic incremental method, in the Markov
randomized incremental method, only one agent updates at any given
time. In particular, an agent in the network updates the current
iterate (by processing locally its own objective function) and either,
passes the new iterate to a randomly chosen neighbor, or, processes it
again.  Thus the order in which the agents update the iterates is
random. This class of incremental algorithms was first proposed in
\cite{Nedic01}, where the agent that receives the iterate is chosen
with uniform probability in each iteration (corresponding to the case
of a fully connected agent network).  Recently, this idea has been
extended in \cite{Johansson07} to the case where the sequence in which
the agents process the information is a time {\it homogeneous Markov
chain}.  The rationale behind this model is that the agent updating
the information at a given time is more likely to pass this
information to a close neighbor rather than to an agent who is further
away.  In this paper, we consider a more general framework than that
of \cite{Johansson07} by allowing the sequence in which the agents
process the information to be a time {\it non-homogeneous Markov
chain}.\footnote{The primary motivation to study such a model are
mobile networks where the network connectivity structure is changing
in time and, thus, the set of the neighbors of an agent is
time-varying.}  We prove the algorithm convergence for diminishing
step-size and establish an error bound for a constant step-size.  This
extends the results in \cite{Johansson07}, where an error bound for a
homogeneous Markov randomized incremental subgradient method is
discussed for a constant step-size and error-free case.

The Markov randomized incremental algorithm is also related to the
decentralized computation model in \cite{Tsitsiklis86,Bertsekas97} for
stochastic optimization problems. However, the emphasis in these
studies is on parallel processing where each agent completely knows
the entire objective function to be minimized.  More closely related
is the work in studies in \cite{Nedic08b} that develops a ``parallel''
version of the unconstrained incremental subgradient algorithm. Also
related is the constrained consensus problem studied in
\cite{Nedic08c} where agents are interested in obtaining a solution to
a feasibility problem, when different parts of the problem are known
to different agents. At a much broader scale, the paper is also
related to the literature on distributed averaging and consensus
algorithms \cite{Tsitsiklis84,Tsitsiklis86,Bertsekas97,
Kar07,Olfati05,Xiao07,Nedic08b,Nedic08c,Olshevsky08,Jadbabaie03}.

Our main contributions in this paper are the development and analysis
of the Markov randomized incremental method with stochastic
subgradients and the use of a time non-homogeneous Markov model for
the sequence of computing agents.  In addition, To the best of our
knowledge, this is among the few attempts made at studying the effects
of stochastic errors on the performance of decentralized optimization
algorithms. The other studies are
\cite{Tsitsiklis84,Tsitsiklis84b,Li87}, \ram{but the algorithms
considered are fundamentally different from the incremental algorithms
studied in this paper}.\footnote{In that work, the components of the
decision vector are distributed while the objective function is known
to all agents. In contrast, in this paper, the objective function data
is distributed, while each agent has an estimate of the entire
decision vector.}

The paper is organized as follows. In Section \ref{sec:problem}, we
formulate the problem of interest, and introduce the cyclic
incremental and Markov randomized incremental method with stochastic
errors.  We also discuss some applications that motivate our interest
in these methods.  In Section \ref{sec:incgrad}, we analyze
convergence properties of the cyclic incremental method. We establish
convergence of the method under diminishing step-size and provide an
error bound for the method with a constant step-size, both valid with
probability 1.  We establish analogous results for the Markov
randomized incremental method in Section \ref{sec:markov}.  We give
some concluding remarks in Section \ref{sec:discuss}.

\section{Problem Formulation and Motivation}
\label{sec:problem}
We consider a network of $m$ agents, indexed by $i=1,\ldots,m$.
The network objective is to solve the following problem:
\begin{equation}\begin{array}{ll}
  \mbox{minimize \  } & f(x)=\sum_{i=1}^{m} f_i(x) \cr
  \mbox{subject to \  } & x \in X,  
 \label{eqn:problem}
\end{array}
\end{equation}
where $x\in\rn$ is a decision or a parameter vector, $X$ is a closed
and convex subset of $\Re^n$, and each $f_i$ is a convex function from
$\rn$ to $\Re$ that is known only to agent $i.$ Problems with the
above structure arise in the context of estimation in sensor networks
\cite{Rabbat05,Sundhar07b}, where $x$ is an unknown parameter to be
estimated and $f_i$ is the cost function that is determined by the
$i$-th sensor's observations (for example, $f_i$ could be the
log-likelihood function in maximum likelihood
estimation). Furthermore, problems with such structure also arise in
resource allocation in data networks. In this context, $x$ is the
resource vector to be allocated among $m$ agents and $f_i$ is the
utility function for agent $i$ \cite{Eryilmaz06}. We discuss these in
more detail later.

To solve the problem (\ref{eqn:problem}) in a network where agents are
connected in a directed ring structure, we consider the cyclic
\emph{incremental subgradient method} \cite{Nedic01}.  Time is slotted
and in each time slot, the estimate is passed by an agent to the next
agent along the ring. In particular, agent $i,$ receives the iterate
from agent $i-1,$ and updates the received estimate using a
subgradient of its ``own objective function $f_i$''.  The updated
iterate is then communicated to the next agent in the cycle, which is
agent $i+1$ when $i<m$ and agent $1$ when $i = m.$ We are interested
in the case where the agent subgradient evaluations have random
errors. Formally, the algorithm is given as follows:
\begin{equation}
\begin{array} {ll}
 &z_{0,k+1} = z_{m,k} = x_{k}, \cr 
 &z_{i,k+1} = \mathcal{P}_{X} \left[ z_{i-1,k+1} - \alpha_{k+1} 
  \left(\nabla
    f_{i}(z_{i-1,k+1})+\e_{i,k+1} \right)\right], 
\end{array}
\label{eqn:incgrad}
\end{equation}
where the initial iterate $x_0 \in X$ is chosen at random.  The vector
$x_{k}$ is the estimate at the end of cycle $k,$ $z_{i,k+1}$ is the
intermediate estimate obtained after agent $i$ updates in $k+1$-st
cycle, $\nabla f_i(x)$ is the subgradient of $f_i$ evaluated at $x,$
and $\e_{i,k+1}$ is a random error. The scalar $\alpha_{k+1}$ is a
positive step-size and $\mathcal{P}_{X}$ denotes Euclidean projection
onto the set $X.$ We study the convergence properties of method
(\ref{eqn:incgrad}) in Section \ref{sec:incgrad} for diminishing and
constant step-size.

In addition, for a network of agents with arbitrary connectivity, we
consider an incremental algorithm where the agent that updates is
selected randomly according to a distribution depending on the agent
that performed the most recent update.  Formally, in this method the
iterates are generated according to the following rule:
\begin{equation}
  x_{k+1} = \mathcal{P}_{X} \left[ x_{k} - \alpha_{k+1} \left(\nabla
  f_{s(k+1)}(x_{k}) +\e_{s(k+1),k+1} \right) \right],
  \label{eqn:markov}
\end{equation}
where the initial iterate $x_0\in X$ is chosen at
random and the agent $s(0)$ that initializes the method is also
selected at random.  The integer $s(k+1)$ is the index of the agent
that performs the update at time $k+1$, and the sequence $\{s(k)\}$ is
modeled as a time non-homogeneous Markov chain with state space
$\{1,\ldots,m\}$. In particular, if agent $i$ was processing at time
$k$, then the agent $j$ will be selected to perform the update at time
$k+1$ with probability $[P(k)]_{i,j}$. Formally, we have
\[\Prob{\{s(k+1)=j\mid s(k)=i\}}=[P(k)]_{i,j}.\]
When there are no errors ($\e_{s(k+1),k+1}=0$) and the
probabilities $[P(k)]_{i,j}$ are all equal to $\frac{1}{m}$,
the method in (\ref{eqn:markov}) coincides with
the incremental method with randomization
that was proposed and studied in \cite{Nedic01}.

Following \cite{Johansson07}, we refer to the method in
(\ref{eqn:markov}) as the {\it Markov randomized incremental
stochastic} algorithm. We analyze convergence properties of this
method in Section \ref{sec:markov} for diminishing and constant
step-sizes.

\subsection{Motivation}
As mentioned, we study the convergence properties of the incremental
algorithms (\ref{eqn:incgrad}) and (\ref{eqn:markov}) for diminishing
and constant step-size, and for zero and non-zero mean errors.  Such
errors may arise directly as computational round-off errors, which are
of interest when the entire network is on a single chip and each agent
is a processor on the chip \cite{Narayanan08}. In addition, stochastic
errors also arise in the following context.

Let the function $f_i(x)$ have the following form
\begin{align*}
f_i(x) = \EXP{g_i(x,R_i)},
\end{align*}
where $\EXP{\cdot}$ denotes the expectation, $R_i\in\Re^d$ is a random
vector and $g_i:\Re^{n\times d}\to\Re$. Agent $i$ does not know the
statistics of $R_i,$ and thus does not know its complete objective
function $f_i$.  However, agent $i$ sequentially observes independent
samples of $R_i$ and uses these samples to determine an approximate
subgradient using the Robbins-Monro approximation \cite{Robbins51} or
Kiefer-Wolfowitz approximation \cite{Kiefer52}. These approximate
sub-gradients can be considered to be the actual sub-gradient
corrupted by stochastic errors.

We next discuss some specific problems that fall within the framework
that we consider and can be solved using the proposed methods.
\subsubsection*{Distributed Regression}
Consider $m$ sensors that sense a time invariant spatial field. Let
$r_{i,k}$ be the measurement made by $i^{th}$ sensor in time slot $k.$
Let $s_i$ be the location of the $i^{th}$ sensor.  Let $h(s;x)$ be a
set of candidate models for the spatial field that are selected based
on a~priori information and parameterized by $x$.  Thus, for each $x,$
the candidate $h(s,x)$ is a model for the spatial field and $h(s_i,x)$
is a model for the measurement $r_{i,k}.$ The problem in regression is
to choose the best model among the set of candidate models based on
the collected measurements $r_{i,k}$, i.e., to determine the value for
$x$ that best describes the spatial field.  In least square
regression, the parameter value $x^*$ corresponding to the best model
satisfies the following relation:
\[
x^*\in \Argmin_{x \in X}\, \lim_{N \to \infty}
\sum_{i=1}^{m}\frac{1}{N} \, \sum_{k=1}^{N} \left( r_{i,k} -
h(s_{i},x)\right)^2.
\]
When the measurements $r_{i,k}$ are corrupted by i.i.d. noise, then 
the preceding relation is equivalent to the following
\[
x^* \in \Argmin_{x \in X} \sum_{i=1}^{m} \EXP{\left(
R_i - h(s_i,x)\right)^2}.
\]
In linear least squares regression, the models $h(s_i,x)$,
$i=1,\ldots,m,$ are linear in $x$, so that each of the functions
$f_i(x)= \EXP{\left( R_i - h(s_i,x)\right)^2}$ is convex in $x$.

\subsubsection*{Distributed Resource Allocation}
An important problem in wireless networks is the fair rate allocation
problem \cite{Eryilmaz06}. Consider a wireless network represented by
a graph with a set of directed (communication) links.  Suppose that
there are $m$ flows indexed $1,\ldots,m,$ whose rate can be adjusted,
and let $x_i$ be the rate of the $i^{th}$ flow.  Each flow is
characterized by a source node $b(i)$ and a destination node $e(i).$
The rate vector $x$ must satisfy some constraints that are imposed by
the individual link capacities of the network. For example, if there
are multiple flows (or, parts of the flows) that use a link of total
capacity $c$ then the total sum of the rates of the flow routed
through that link must not be greater than $c.$ Apart from this there
could also be some queuing delay constraints. Thus, only flow vectors
that are constrained to a set $X$ can be routed through the
network. Associated with each flow, there is a reward function
$U_i(x_i)$ depending only on the flow rate $x_i$ and known only at the
source node $b(i).$ \ram{The reward function is typically a concave
and increasing function.} In the \emph{fair rate allocation problem},
the source nodes $\{b(\ell)\}$ need to determine the optimal flow rate
$x^*$ that maximizes the total network utility. Mathematically, the
problem is to determine $x^*$ such that
\[
x^* \in \Argmax_{x \in X} \sum_{i=1}^{m} U_{i}(x_i).
\]
In some networks, the same flow can communicate different types of
traffic that has different importance in different time slots.  For
example, in an intruder detection network, a ``detected'' message is
more important (and is rewarded/weighted more) than a ``not detected''
message or some other system message.  Thus, the reward function is
also a function of the contents of the flow: if the type of flow $i$
in time slot $k$ is $r_{i,k},$ where $r_{i,k}$ takes values from the
set of all possible types of flow data, then the reward is
$U_i(x_i,r_{i,k})$ at time $k$.  When the type of traffic on each flow
across slots is i.i.d, the \emph{fair allocation rate} problem can be
written as
\[
\max_{x \in X} \sum_{i=1}^{m} \EXP{U_{i}( x_{i}, R_{i})}.
\]
The statistics of $R_i$ may not be known since they may depend upon
external factors such as the frequency of intruders in an intruder
detection network.

\subsection{Notation and Basics}
\label{sec:basic}
We view vectors as columns.  We write $x^Ty$ to denote the inner
product of two vectors $x$ and $y$.  We use $\|\cdot\|$ to denote the
standard Euclidean norm.  For a vector $x$, we use $x_i$ to denote its
$i$-th component.  For a matrix $A,$ we use $[A]_{i,j}$ to denote its
$(i,j)$-th entry, $[A]_i$ its $i$-th row and $[A]^j$ its $j$-th
column. We use $e$ to denote a vector with each entry equal to 1.

We use $f^*$ to denote the optimal value of the problem
(\ref{eqn:problem}), and we use $X^*$ to denote its optimal
set. Throughout the paper, we assume that the optimal value $f^*$ is
finite.

In our analysis, we use the subgradient defining
property. Specifically, for a convex function $f:\rn\to\Re$, the
vector $\nabla f(x)$ is a {\it subgradient of $f$ at $x$ } when the
following relation is satisfied:
\begin{eqnarray}
\nabla f(x)^T (y - x)\le f(y) - f(x) 
\qquad\hbox{for all $y\in \rn$}
\label{eqn:lip3}
\end{eqnarray}
(see, for e.g., \cite{Bertsekas03}).

\section{Cyclic Incremental Subgradient Algorithm}
\label{sec:incgrad}
Recall, that the cyclic incremental stochastic subgradient algorithm
is given by
\begin{equation}
  \begin{array}{ll}
    z_{0,k+1} &= z_{m,k} = x_{k}, \cr
    z_{i,k+1} &=
    \mathcal{P}_{X} \left[ z_{i-1,k+1} - \alpha_{k+1} \left( \nabla
    f_{i}(z_{i-1,k+1}) + \epsilon_{i,k+1} \right)
    \right],
  \end{array}
\label{eqn:stocincgrad}
\end{equation}
where $x_0\in X$ is a random initial vector, $\nabla f_i(x)$ 
is a subgradient of $f_i$ at $x$, $\epsilon_{i,k+1}$ is a
random noise vector and $\alpha_{k+1}>0$ is a step-size.

The main difficulty in the study of the incremental stochastic
subgradient algorithm is that the expected direction in which the
iterate is adjusted in each sub-iteration is not necessarily a
subgradient of the objective function $f$. For this reason, we cannot
directly apply the classic stochastic approximation convergence
results of \cite{Kushner78,Polyak87,Ermoliev76} to study the
convergence of method in (\ref{eqn:incgrad}).  The key relation in our
analysis is provided in Lemma \ref{lemma:key} in Section
\ref{subsec:prelim}.  Using this lemma and a standard super-martingale
convergence result, we obtain results for diminishing step-size in
Theorem \ref {thm:QISGA1}.  Furthermore, by considering a related
``stopped'' process to which we apply a standard supermartingale
convergence result, we obtain the error bound results for a constant
step-size in Theorems \ref{thm:igac1} and \ref{thm:igac2}.

We make the following basic assumptions on the set $X$ and the functions $f_i$.
\begin{assumption}\label{ass:convex} 
The set $X\subseteq\rn$ is closed and convex.
The function $f_i:\rn \to\Re$ is convex for each $i\in\{1,\ldots,m\}.$
\end{assumption}

In our analysis, we assume that the first and the second moments of
the subgradient noise $\e_{i,k}$ are bounded uniformly over the
agents, conditionally on the past realizations.  In particular, we
define $F_{k}^{i}$ as the $\sigma$-algebra generated by $x_0$ and the
subgradient errors $\epsilon_{1,1}, \ldots, \epsilon_{i,k}$, and
assume the following.
\begin{assumption}
  \label{ass:erindep}
There exist deterministic scalar sequences 
$\{\mu_k\}$ and $\{\nu_{k}\}$ such that 
\[
\|\EXP{\epsilon_{i,k}\mid F_{k}^{i-1}}\|
\le\mu_k\qquad\hbox{for all $i$ and $k$},\] 
\[\EXP{\|\epsilon_{i,k}\|^2 \mid F_{k}^{i-1} } 
\le \nu^2_{k}
\qquad\hbox{for all $i$ and $k$}.\]
\end{assumption}

Assumption \ref{ass:erindep} holds for example, when the errors
$\epsilon_{i,k}$ are independent across both $i$ and $k,$ and have
finite moments. Note that under the assumption that the second moments
are bounded, from Jensen's inequality we readily have
\begin{equation}
\left\| \EXP{\epsilon_{i,k} \mid F_{k}^{i-1}} \right\| \leq
    \sqrt{\EXP{\|\epsilon_{i,k}\|^2 \mid F_{k}^{i-1} }} \leq \nu_{k}.
   \label{eqn:jensen}
\end{equation}
However, for a constant step-size, the terms $\left\|
\EXP{\epsilon_{i,k} \mid F_{k}^{i-1}} \right\|$ and
$\EXP{\|\epsilon_{i,k}\|^2 \mid F_{k}^{i-1} }$ affect the error bounds
on the performance of the incremental method differently, as seen in
Section~\ref{ssec:css1}.  For this reason, we prefer to use different
upper-bounds for the terms $\left\| \EXP{\epsilon_{i,k} \mid
F_{k}^{i-1}} \right\|$ and $\EXP{\|\epsilon_{i,k}\|^2 \mid F_{k}^{i-1}
}.$ We will also, without any loss of generality, assume that $\mu_k <
\nu_k.$

We also assume that the subgradients $\nabla f_i(x)$ are uniformly
bounded over the set $X$ for each $i$. This assumption is commonly
used in the convergence analysis of subgradient methods with a
diminishing or a constant step-size.
\begin{assumption}\label{ass:sgdbound}
For every $i$, the subgradient set of the function $f_i$ at $x\in X$
is nonempty and uniformly bounded over the set $X$ by a constant
$C_i$, i.e.,
$$\|\nabla f_i(x)\|\le C_i\qquad
\hbox{for all subgradients $\nabla f_i(x)$}\quad
\hbox{and for all $x\in X$}.$$
\end{assumption}

Assumption \ref{ass:sgdbound} holds for example, when each $f_i$ is a
polyhedral function or when the set $X$ is compact. 

\subsection{Preliminaries}
\label{subsec:prelim}
In this section, we provide a lemma establishing a basic relation for
the iterates generated by the incremental method
(\ref{eqn:stocincgrad}) and any step-size rule. This relation plays a
key role in our subsequent development.

\begin{lemma}
  \label{lemma:key}
  Let Assumptions~\ref{ass:convex}, \ref{ass:erindep}, and
  \ref{ass:sgdbound} hold.  Then, the iterates generated by algorithm
  (\ref{eqn:stocincgrad}) are such that for any step-size rule and for
  any $y\in X$, 
    \begin{align}
      \EXP{\| d_{k+1} (y) \|^2 \mid F_{k}^{m} } \leq& \|d_k (y) \|^2 -
      2 \alpha_{k+1} \left( f(x_k) - f(y) \right) \nonumber \\ 
      & +
      2\alpha_{k+1} \mu_{k+1}\sum_{i=1}^m \EXP{\| d_{i-1,k+1}(y)
      \| \mid F_{k}^{m} } \nonumber \\ 
      & + \alpha_{k+1}^2 \left(m
      \nu_{k+1} + \sum_{i=1}^{m} C_i\right)^2,
      \label{lem}
    \end{align}
    where $d_k(y) = x_k - y$ and $d_{i,k+1}(y) = z_{i,k+1} - y$
    for all $k$.
\end{lemma}

\begin{proof}
Using the iterate update rule in (\ref{eqn:stocincgrad}) and the 
non-expansive property of the
  Euclidean projection,
  we obtain for any $y\in X$,
  \begin{align*}
    \|d_{i,k+1}(y)\|^2 =& \left\|\mathcal{P}_{X} \left[ z_{i-1,k+1} -
      \alpha_{k+1} \nabla f_i(z_{i-1,k+1}) - \alpha_{k+1}
      \epsilon_{i,k+1} \right] - y \right\|^2 \nonumber \\ \leq&
      \left\| z_{i-1,k+1} - \alpha_{k+1} \nabla f_i(z_{i-1,k+1}) -
      \alpha_{k+1} \epsilon_{i,k+1} - y \right\|^2 \nonumber \\ =&
      \|d_{i-1,k+1}(y)\|^2 - 2 \alpha_{k+1} d_{i-1,k+1}(y)^T \nabla
      f_{i}(z_{i-1,k+1}) \\& - 2 \alpha_{k+1} d_{i-1,k+1}(y)^T
      \epsilon_{i,k+1} + \alpha_{k+1}^2\left\| \epsilon_{i,k} +
      \nabla f_{i}(z_{i-1,k+1}) \right\|^2.
  \end{align*}
  Taking conditional expectations with respect to the $\s$-field
  $F_{k+1}^{i-1}$, we further obtain 
  \begin{align}
    \EXP{\|d_{i,k+1}(y)\|^2 \mid F_{k+1}^{i-1} } \leq&
    \|d_{i-1,k+1}(y)\|^2 - 2 \alpha_{k+1} d_{i-1,k+1}(y)^T \nabla f_i(
    z_{i-1,k+1}) \nonumber \\ &- 2 \alpha_{k+1} d_{i-1,k+1}(y)^T
    \EXP{\epsilon_{i,k+1} \mid F_{k+1}^{i-1}} \nonumber \\
    &+ \alpha_{k+1}^2\EXP{\left\| \epsilon_{i,k+1} +\nabla
    f_{i}(z_{i-1,k+1}) \right\|^2 \mid F_{k+1}^{i-1}}. 
    \label{eqn:step1}
    \end{align}
    
  We now estimate the last two terms in the right hand side of 
  the preceding equation by using Assumption~\ref{ass:erindep} 
  on the error moments. 
  In particular, we have for all $i$ and $k$,
  {\small
  \[- d_{i-1,k+1}(y)^T
    \EXP{\epsilon_{i,k+1} \mid F_{k+1}^{i-1}} 
   \le \|d_{i-1,k+1}(y)\|\,
    \left\| \EXP{\epsilon_{i,k+1} \mid F_{k+1}^{i-1}} \right\|
   \le \m_{k+1} 
    \|d_{i-1,k+1}(y)\|.\] 
   }

  Next, we estimate the last term in Eq.\ (\ref{eqn:step1})
  by using Assumption~\ref{ass:erindep} on the error moments in
  and Assumption~\ref{ass:sgdbound} on the subgradient norms.  
  We have all $i$ and $k$,
  \begin{align*}
      \EXP{ \left\| \epsilon_{i,k+1} + \nabla
      f_{i}(z_{i-1,k+1}) \right\|^2 \mid F_{k+1}^{i-1}} =& \EXP{
      \left\| \epsilon_{i,k+1} \right\|^2 \mid F_{k+1}^{i-1}} +
      \left\|\nabla f_{i}(z_{i-1,k+1}) \right\|^2 \\&+ 2 \nabla
      f_{i}(z_{i-1,k+1})^T \EXP{ \epsilon_{i,k+1} \mid F_{k+1}^{i-1}}
      \\ =& \EXP{ \left\| \epsilon_{i,k+1} \right\|^2 \mid
      F_{k+1}^{i-1}} + \left\|\nabla f_{i}(z_{i-1,k+1}) \right\|^2
      \\&+ 2 \|\nabla f_{i}(z_{i-1,k+1})\| \|\EXP{ \epsilon_{i,k+1}
      \mid F_{k+1}^{i-1}}\| \\\leq& ( \nu_{k+1} + C_i )^2,
  \end{align*}
  where in the last inequality we use
  $\EXP{ \left\| \epsilon_{i,k} \right\|^2 \mid
      F_{k}^{i-1}}\le \nu_{k}^2$ and 
  $\left\| \EXP{\epsilon_{i,k} \mid F_{k}^{i-1}} \right\| \le \nu_{k}$
  for all $k$
  [cf.\ Eq.\ (\ref{eqn:jensen})].
  Combining the preceding two relations and the inequality in 
  (\ref{eqn:step1}), 
  we obtain for all $y\in X$,
  \begin{align}
      \EXP{\|d_{i,k+1}(y)\|^2 \mid F_{k+1}^{i-1} } \leq&
      \|d_{i-1,k+1}(y)\|^2 - 2 \alpha_{k+1} d_{i-1,k+1}(y)^T \nabla
      f_i( z_{i-1,k+1}) \nonumber \\& + 2 \alpha_{k+1} \mu_{k+1} \|
      d_{i-1,k+1}(y) \| + \alpha_{k+1}^2(\nu_{k+1} + C_i)^2. 
      \label{eqn:ineq2}
  \end{align} 
  
  We now estimate the second term in the right hand side of 
  the preceding relation. From the subgradient inequality
  in (\ref{eqn:lip3}) we have 
  {\small
  \begin{eqnarray}
    - d_{i-1,k+1}(y)^T \nabla f_i( z_{i-1,k+1}) &=& - (z_{i-1,k+1} - y)^T
      \nabla f_i( z_{i-1,k+1}) \cr
    &\leq& -\left(
    f_i(z_{i-1,k+1}) - f_i(y) \right) \cr
    &=& - \left(
    f_i(x_{k}) - f_i(y) \right) - \left( f_i(z_{i-1,k+1}) - f_i(x_k)
    \right) \cr
    &\leq& -\left( f_i(x_{k}) - f_i(y) \right) -
    \left( \nabla f_i(x_k) \right)^T \left( z_{i-1,k+1} - x_k \right)
    \label{eqn:smart} \\
    &\leq& - \left( f_i(x_{k}) - f_i(y) \right) +
    C_i \left\| z_{i-1,k+1} - x_k \right\|. 
    \label{eqn:smart2}
    \end{eqnarray}} 
    In (\ref{eqn:smart}) we have again used the subgradient
    inequality (\ref{eqn:lip3}) to bound $f_i(z_{i-1,k+1}) -
    f_i(x_k),$ while in (\ref{eqn:smart2}) we have used 
    the subgradient norm bound from Assumption~~\ref{ass:sgdbound}.  
    We next consider the term $\left\| z_{i-1,k+1} - x_k \right\|$.
    From (\ref{eqn:stocincgrad}) we have
    \[ \| z_{i-1,k+1} - x_k \| 
    = \left\|\sum_{j=1}^{i-1} \left(z_{j,k+1} - z_{j-1,k+1}\right) \right\|
    \leq \sum_{j=1}^{i-1} \| z_{j,k+1} -
    z_{j-1,k+1}\|.\]
    By the non-expansive property of the projection, we further have
    \begin{align}
    \| z_{i-1,k+1} - x_k \| \leq \alpha_{k+1} \sum_{j=1}^{i-1}
    \left(\|\nabla f_{j}(z_{j-1,k+1}) \| + \| \epsilon_{j,k+1}
    \|\right) \leq \alpha_{k+1} \sum_{j=1}^{i-1} \left(C_j +
    \|\epsilon_{j,k+1} \|\right).\label{eqn:k3}
     \end{align}
    By combining the preceding relation with Eq.\ (\ref{eqn:smart2}),
    we have
    \[
    - d_{i-1,k+1}(y)^T \nabla f_i( z_{i-1,k+1}) \leq - \left(
    f_i(x_{k}) - f_i(y) \right) + \alpha_{k+1} C_i \sum_{j=1}^{i-1}
    \left(C_j + \|\epsilon_{j,k+1} \| \right).
    \]
    By substituting the preceding estimate in the inequality in
    (\ref{eqn:ineq2}), we obtain for all $y\in X$,
    \begin{align*}
      \EXP{\|d_{i,k+1}(y)\|^2 \mid F_{k+1}^{i-1} } \leq&
      \|d_{i-1,k+1}(y)\|^2 - 2 \alpha_{k+1} \left(f_i(x_k) - f_i(y)
      \right)\\ &+2\a_{k+1}^2C_i\sum_{j=1}^{i-1} \left(C_j +
      \|\epsilon_{j,k+1} \| \right)\\ &+ 2
      \alpha_{k+1} \mu_{k+1} \|d_{i-1,k+1}(y) \| + \alpha^2_{k+1} (
      C_i + \nu_{k+1} )^2 .
    \end{align*}
    Taking the expectation conditional on $F_{k}^m$, we obtain
    \begin{align*}
      \EXP{\|d_{i,k+1}(y)\|^2 \mid F_k^m } \leq&
      \EXP{\|d_{i-1,k+1}(y)\|^2\mid F_k^m} - 2 \alpha_{k+1}
      \left(f_i(x_k) - f_i(y) \right)\\ &+ 2 \alpha_{k+1} \mu_{k+1}
      \EXP{\|d_{i-1,k+1}(y) \| \mid F_k^m } \\
      &+2\a_{k+1}^2C_i\sum_{j=1}^{i-1} \left(C_j + \nu_{k+1}\right) +
      \alpha^2_{k+1} ( C_i + \nu_{k+1} )^2,
  \end{align*}
  where we have used Assumption~2 and Jensen's inequality to bound
  $\EXP{\|\epsilon_{j,k+1} \| \mid F_k^m}$ by $\nu_{k+1}$ [cf.\ Eq.\
  (\ref{eqn:jensen})].  Summing over $i=1,\ldots,m,$ and noting that
  $d_{0,k+1} (y)=x_k-y$, we see that
  \begin{align*}
    \EXP{\|d_{k+1}(y)\|^2 \mid F_k^m } \leq& \|d_k(y)\|^2 - 2
      \alpha_{k+1} \left(f(x_k) - f(y) \right)\\ &+ 2 \alpha_{k+1}
      \mu_{k+1} \sum_{i=1}^m 
     \EXP{\|d_{i-1,k+1}(y) \| \mid F_k^m
      }\\ &+2\a_{k+1}^2\sum_{i=1}^m C_i\sum_{j=1}^{i-1} \left(C_j +
      \nu_{k+1}\right) + \sum_{i=1}^m\alpha^2_{k+1} ( C_i + \nu_{k+1}
      )^2.
  \end{align*}
  Finally, by noting that 
  \[ 2\sum_{i=1}^m 
  C_i\sum_{j=1}^{i-1} \left(C_j + \nu_{k+1}\right) + \sum_{i=1}^m
  ( C_i + \nu_{k+1} )^2=\left( \sum_{i=1}^{m} C_i+ m \nu_{k+1}
  \right)^2,\] we obtain for all $y\in X$, and all $i$ and $k$,
  \begin{align*}
    \EXP{\|d_{k+1}(y)\|^2 \mid F_{k}^{m} } \leq& \|d_{k}(y)\|^2 - 2
    \alpha_{k+1} \left(f(x_k) - f(y) \right) \\& + 2 \alpha_{k+1}
    \mu_{k+1} \sum_{i=1}^m 
    \EXP{\|d_{i-1,k+1}(y) \| \mid F_k^m }
    \\& + \alpha^2_{k+1} \left( \sum_{i=1}^{m} C_i+ m \nu_{k+1}
    \right)^2.
  \end{align*}
\end{proof}

\subsection{Convergence for diminishing step-size}
We here study the convergence of the method in (\ref{eqn:stocincgrad})
for diminishing step-size rule. In our analysis, 
we use the following
result due to Robbins and Siegmund (see Lemma~11, Chapter~2.2,
\cite{Polyak87}).
\begin{lemma}
\label{lemma:asm}
 Let $(\Omega, \mathcal{F},\mathcal{P})$ be a probability space and let
 $\mathcal{F}_0 \subset \mathcal{F}_1 \subset \ldots$ be a sequence of
 sub $\sigma$-fields of $\mathcal{F}.$ Let $u_k, v_k$ and $w_k,$
 $k=0,1,2\ldots,$ be non-negative $\mathcal{F}_k$-measurable random
 variables and let $\{q_k\}$ be a deterministic sequence. Assume that
 $\sum_{k=0}^{\infty} q_k <\infty,$ and $\sum_{k=0} w_k <\infty$ and
\[
\EXP{u_{k+1} \mid {\cal F}_k } \leq (1 + q_{k}) u_{k} - v_{k} +w_{k}
\]
hold with probability~1.  Then, with probability 1, the sequence
$\{u_{k}\}$ converges to a non-negative random variable and
$\sum_{k=0}^{\infty} v_{k} < \infty$.
\end{lemma}

We next provide a convergence result for diminishing step-sizes. 
\begin{theorem}
  \label{thm:QISGA1}
  Let Assumptions~\ref{ass:convex},
  \ref{ass:erindep} and \ref{ass:sgdbound} hold.  
  Assume that the step-size sequence
  $\{\a_k\}$ is positive and such that $\sum_{k=1}^{\infty} \alpha_k =
  \infty$ and $\sum_{k=1}^{\infty} \alpha_k^2 < \infty$. In addition,
  assume that the bounds $\mu_{k}$ and $\nu_{k}$ on the moments of the
  error sequence $\{\e_{i,k}\},$ are such that
  \[\sum_{k=1}^{\infty} \alpha_k\mu_k < \infty,\qquad
    \sum_{k=1}^{\infty} \alpha_k^2\nu_k^2 < \infty.
  \]
  Also, assume that the optimal set $X^*$ is nonempty.  Then, the
  iterate sequence $\{x_k\}$ generated by the method
  (\ref{eqn:stocincgrad}) converges to an optimal solution with
  probability~1.
\end{theorem}
\begin{proof}
First note that all the assumptions of Lemma~\ref{lemma:key} are
satisfied.  Let $x^*$ be an arbitrary point in $X^*.$ By letting $y =
x^*$ in Lemma~\ref{lemma:key}, we obtain for any $x^*\in X^*$,
\begin{align}
  \EXP{\|d_{k+1}(x^*)\|^2 \mid F_{k}^{m}} \leq& \|d_k (x^*) \|^2 - 2
  \alpha_{k+1} \left(f(x_k) - f^* \right) \nonumber\\&+ 2\alpha_{k+1}
  \mu_{k+1} \sum_{i=1}^m \EXP{\|d_{i-1,k+1}(x^*) \| \mid F_k^m }
  \nonumber\\&+ \alpha_{k+1}^2 \left(m \nu_{k+1} + \sum_{i=1}^{m}
  C_i\right)^2. \label{eqn:k4}
\end{align}
We relate $\|d_{i-1,k+1} (x^*) \|$ to $\|d_k (x^*)\|$ by using the
triangle inequality of norms,
\begin{align*}
\|d_{i-1,k+1} (x^*) \| = \| z_{i-1,k+1} - x_{k} + x_{k} - x^* \|
\leq \| z_{i-1,k+1} - x_{k}\| + \|d_k (x^*)\|.
\end{align*}
Substituting for $\| z_{i-1,k+1} - x_{k}\|$ from (\ref{eqn:k3}) we
obtain
\begin{align*}
\|d_{i-1,k+1} (x^*) \| &\leq \alpha_{k+1} \sum_{j=1}^{i-1} \left(C_j +
    \|\epsilon_{j,k+1}\|\right) + \|d_k (x^*)\|.
\end{align*}
Taking conditional expectations, we further obtain
\begin{align*}
\EXP{\|d_{i-1,k+1} (x^*) \| \mid F_{k}^{m}} \leq& \|d_k (x^*)\| +
\alpha_{k+1} \sum_{j=1}^{i-1} \left(C_j + \nu_{k+1}\right),
\end{align*}
where we have used Assumption~\ref{ass:erindep} and Jensen's
inequality to bound $\EXP{\|\epsilon_{j,k+1}\| \mid F_{k}^{m}}$ by
$\nu_{k+1}.$ Using the preceding inequality in (\ref{eqn:k4}), we have
\begin{align*}
  \EXP{\|d_{k+1}(x^*)\|^2 \mid F_{k}^{m}} 
  \leq& 
  \|d_k (x^*) \|^2 - 2
  \alpha_{k+1} \left(f(x_k) - f^* \right) \\& + 2m \alpha_{k+1} \mu_{k+1}
  \|d_{k}(x^*) \|  + 2 \alpha^2_{k+1} \mu_{k+1} \sum_{i=1}^{m}
  \sum_{j=1}^{i-1} \left(C_j + \nu_{k+1}\right) \\& + \alpha_{k+1}^2
  \left(m \nu_{k+1} + \sum_{i=1}^{m} C_i\right)^2.
\end{align*}
Next, using the inequality
\[
2  \|d_{k}(x^*) \| \leq 1 +  \|d_{k}(x^*) \|^2,
\]
we obtain
\begin{align}
  \EXP{\|d_{k+1}(x^*)\|^2 \mid F_{k}^{m}} \leq& \l(1 + m
  \alpha_{k+1}\mu_{k+1} \right) \|d_k (x^*) \|^2 - 2 \alpha_{k+1}
  \left(f(x_k) - f^* \right) \nonumber \\& + m \alpha_{k+1} \mu_{k+1}
  + 2 \alpha^2_{k+1} \mu_{k+1} \sum_{i=1}^{m} \sum_{j=1}^{i-1}
  \left(C_j + \nu_{k+1}\right) \nonumber \\& + \alpha_{k+1}^2 \left(m
  \nu_{k+1} + \sum_{i=1}^{m} C_i\right)^2. \label{eqn:newlem}
\end{align}
By the assumptions on the step-size, and the sequences $\{\mu_k\}$ and
$\{\nu_k\}$, we further have
\begin{align*}
&\sum_{k=0}^{\infty} m \alpha_{k+1} \mu_{k+1} < \infty, \\
&\sum_{k=0}^{\infty} 2 \alpha^2_{k+1} \mu_{k+1} \sum_{i=1}^{m}
\sum_{j=1}^{i-1} \left(C_j + \nu_{k+1}\right) \leq 2
\sum_{k=0}^{\infty} \sum_{i=1}^{m} \sum_{j=1}^{i-1} \l( \alpha^2_{k+1}
\mu_{k+1} C_j + \alpha^2_{k+1} \nu^2_{k+1} \r) < \infty, \\
&\sum_{k=0}^{\infty} \alpha_{k+1}^2 \left(m \nu_{k+1} + \sum_{i=1}^{m}
C_i\right)^2 \le 2 \sum_{k=0}^{\infty} \alpha_{k+1}^2
\left(m^2\nu^2_{k+1} + \left(\sum_{i=1}^{m}
C_i\right)^2\right)<\infty,
\end{align*}
where in the second relation above, we have used
$\mu_{k+1}\le\nu_{k+1}$ [cf.\ Eq.\ (\ref{eqn:jensen})], while in the
last inequality, we have used $(a+b)^2\le 2 (a^2 + b^2)$ valid for any
scalars $a$ and~$b$.  Thus, the conditions of Lemma~\ref{lemma:asm}
are satisfied with $u_{k}= \|d_k (x^*) \|^2,$ $\mathcal{F}_k =
F^{m}_{k},$ $q_{k} = m \alpha_{k+1} \mu_{k+1},$ $v_{k} = 2
\alpha_{k+1} \left(f(x_k) - f^* \right)$ and
\begin{align*}
  w_{k} = m \alpha_{k+1} \mu_{k+1} + 2 \alpha^2_{k+1} \mu_{k+1}
  \sum_{i=1}^{m} \sum_{j=1}^{i-1} \left(C_j + \nu_{k+1}\right) +
  \alpha_{k+1}^2 \left(m \nu_{k+1} + \sum_{i=1}^{m} C_i\right)^2.
\end{align*}
Therefore, with probability 1, the scalar $\|d_{k+1}(x^*)\|^2$
converges to some non-negative random variable for every $x^*\in
X^*$. Also with probability 1, we have
\[
\sum_{k=0}^{\infty} \alpha_{k+1} \left(f(x_k) - f^* \right) < \infty.
\]
Since $\sum_{k=1}^{\infty} \alpha_k = \infty,$ it follows that
$\liminf_{k\to\infty}f(x_k) = f^*$ with probability 1.  By considering
a sample path for which $\liminf_{k\to\infty}f(x_k) = f^*$ and
$\|d_{k+1}(x^*)\|^2$ converges for any $x^*$, we conclude that the
sample sequence must converge to some $x^*\in X^*$ in view of
continuity of $f$.  Hence, the sequence $\{x_k\}$ converges to some
vector in $X^*$ with probability 1.
\end{proof}

Note that under assumptions of Theorem \ref{thm:QISGA1}, it can be
seen that $\EXP{\di{x_k}^2}$ also converges to $0.$ In particular,
since the solution set $X^*$ is closed and convex, there exists a
point $x_k^* \in X^*$ that is closest to $x_k$ for every $k$.  Letting
$y=x_k^*$ in relation (\ref{eqn:newlem}) and using the fact that
$\di{x_{k+1}} \leq d_{k+1}(x_k^*)$ with probability $1$, we obtain
for all $k$,
\begin{align*}
  \EXP{\di{x_{k+1}}^2 \mid F_k^m} \leq& \l(1 + m
  \alpha_{k+1}\mu_{k+1} \right)\di{x_k}^2 - 2 \alpha_{k+1}
  \left(f(x_k) - f^* \right) \nonumber \\& + m \alpha_{k+1} \mu_{k+1}
  + 2 \alpha^2_{k+1} \mu_{k+1} \sum_{i=1}^{m} \sum_{j=1}^{i-1}
  \left(C_j + \nu_{k+1}\right) \nonumber \\& + \alpha_{k+1}^2 \left(m
  \nu_{k+1} + \sum_{i=1}^{m} C_i\right)^2.
\end{align*}
Taking expectations, we obtain for all $k$,
\begin{align*}
  \EXP{\di{x_{k+1}}^2} \leq& \l(1 + m \alpha_{k+1}\mu_{k+1}
  \right)\EXP{\di{x_k}^2} - 2 \alpha_{k+1} \left(f(x_k) - f^* \right)
  \nonumber \\& + m \alpha_{k+1} \mu_{k+1}+ 2 \alpha^2_{k+1} \mu_{k+1}
  \sum_{i=1}^{m} \sum_{j=1}^{i-1} \left(C_j + \nu_{k+1}\right)
  \nonumber \\& + \alpha_{k+1}^2 \left(m \nu_{k+1} + \sum_{i=1}^{m}
  C_i\right)^2.
\end{align*}
From the deterministic analog of Lemma~\ref{lemma:asm}, we can argue
that $\EXP{\di{x_{k+1}}^2}$ converges and $\liminf_{k \to \infty}
\EXP{f(x_k)} = f^*$.  Since $\{x_k\}$ converges to a point in $X^*$
with probability 1, it follows that $\EXP{\di{x_{k+1}}^2}$ converges
to $0.$

\subsection{Error bound for constant step-size}
\label{ssec:css1}
Here, we study the behavior of the iterates $\{x_k\}$ generated by the
method (\ref{eqn:stocincgrad}) with a constant step-size rule, i.e.,
$\alpha_k = \alpha$ for all $k$.  In this case, we cannot guarantee
the convergence of the iterates, however, we can provide bounds on the
performance of the algorithm. In the following lemma, we provide an
error bound for the expected values $\EXP{f(x_k)}$ and a bound for
$\inf_k f(x_k)$ that holds with probability 1. The proofs of these
results are similar to those used in \cite{Nedic01b}.
 
\begin{theorem}
  \label{thm:igac1}
  Let Assumptions~\ref{ass:convex} and \ref{ass:erindep} hold. Let the
  sequence $\{x_k\}$ be generated by the method
  (\ref{eqn:stocincgrad}) with a constant step-size rule, i.e.,
  $\alpha_k = \alpha$ for all $k\ge1.$  Also, assume that
  the set $X$ is bounded, and
  \[\mu=\sup_{k \ge1}\mu_{k}<\infty,\qquad 
  \nu= \sup_{k \ge 1} 
  \nu_{k} <\infty.\] We then have
  \begin{equation}
    \liminf_{k\to\infty} \EXP{f(x_k)} \leq f^* + m \mu\,\max_{x,y\in X}\|x-y\|
    + \frac{\alpha}{2} \left(\sum_{i=1}^{m} C_i + m\nu\right)^2, 
    \label{eqn:result1b}
  \end{equation}
  and with probability 1,
  \begin{equation}
    \inf_{k\geq 0} f(x_k) \leq f^* + m \mu\,\max_{x,y\in X}\|x-y\| 
    + \frac{\alpha}{2}
     \left(\sum_{i=1}^{m} C_i + m \nu \right)^2.
     \label{eqn:result2b} 
  \end{equation}
\end{theorem}
\begin{proof}
Since $X$ is compact and each $f_i$ is convex over $\rn$,
the subgradients of $f_i$ are bounded over $X$ for each $i$. Thus, all
the assumptions of Lemma \ref{lemma:key} are satisfied. 
Furthermore, the optimal set $X^*$ is non-empty.
Since $\mu_k\le \mu$ and $\nu_k\le \nu$ for 
all $k$, and $\| d_{i-1,k+1}(y) \|\le \max_{x,y\in X}\|x-y\|,$ 
according to the relation of
Lemma~\ref{lemma:key}, 
we have for $y=x^*\in X^*$,
\begin{align}
    \EXP{\|d_{k+1}(x^*)\|^2 \mid F_{k}^{m} } \leq& \|d_{k}(x^*)\|^2 
    - 2\a\left(f(x_k) - f^* \right)
    + 2 m\a\mu \max_{x,y}\|x-y\| \cr
    & + \a^2\left(
    \sum_{i=1}^{m} C_i+ m \nu \right)^2.
    \label{eqn:rel1}
  \end{align}
By taking the total expectation, we obtain for all $y\in X$ and all $k$,
\begin{align*}
  \EXP{\|d_{k+1}(x^*)\|^2} \leq& \EXP{\|d_k (x^*) \|^2} - 2 \alpha
  \left(\EXP{f(x_k)} - f^* \right) + 2m\a\mu \max_{x,y}\|x-y\| \cr
   &+
  \alpha^2 \left( \sum_{i=1}^{m} C_i+ m \nu \right).
\end{align*}
Now, assume that the relation (\ref{eqn:result1b}) does not hold.  
Then there will exist
a $\gamma > 0$ and an index $k_{\gamma}$ such that for all $k
> k_{\gamma}$,
\[
\EXP{f(x_k)} \geq f^* + \gamma  + m\mu\max_{x,y\in X}\|x-y\| +
\frac{\alpha}{2} \left(\sum_{i=1}^{m} C_i + m\nu \right)^2.
\]
Therefore, for $k > k_{\gamma}$, we have
\begin{align*}
  \EXP{\|d_{k+1} (x^*) \|^2 } \leq& \EXP{\|d_k (x^*) \|^2} - 2
  \alpha \left( \gamma + m\mu\max_{x,y\in X}\|x-y\|
   + \frac{\alpha}{2}
  \left(\sum_{i=1}^{m} C_i + m\nu \right)^2\right) \\& + 2m
  \alpha \mu\max_{x,y\in X}\|x-y\|
  + \alpha^2 \left( \sum_{i=1}^{m} C_i+ m \nu\right)^2 \\
  \leq& \EXP{\|d_{k} (x^*) \|^2 } - 2 \alpha \gamma.
\end{align*}
Hence, for $k \geq k_{\gamma},$
\[
  \EXP{\|d_{k+1} (x^*) \|^2 } \leq
  \EXP{\|d_{k_\gamma}(x^*)\|^2 } - 2 \gamma \alpha (k - k_{\gamma}).
\]
For sufficiently large $k,$ the right hand side of the preceding
relation is negative, yielding a contradiction.  Thus the relation
(\ref{eqn:result1b}) must hold.

We now prove the relation in (\ref{eqn:result2b}).  Define the set
\begin{align*}
L_{N} &= \left\{x \in X: f(x) < f^* 
       + \frac{1}{N} + m\mu \max_{x,y\in X}\|x-y\|+
	\frac{\alpha}{2} \left(\sum_{i=1}^{m} C_i + m\nu \right)^2
	\right\}.
\end{align*}
Let $x^* \in X^*$ and 
define the sequence $\hat{x}_{k}$ as follows:
\[
\hat{x}_{k+1} = \begin{cases} x_{k+1} &\mbox{if  $\hat{x}_{k}
  \notin L_N$,} \\ x^* &\mbox{if $\hat{x}_{k} \in L_N$}. \end{cases}
\]
Thus, the process $\{\hat{x}_k\}$ is identical to the process $\{x_k\},$
until $\{x_k\}$ enters the set $L_N.$ Define
\[
\hat{d}_k(y) = \hat{x}_k - y.
\]
Let us first consider the case when $\hat{x}_k \in L_{N}.$ Since 
$\hat x_k=x^*$ and $\hat x_{k+1}=x^*$, we have $ \hat{d}_{k}(x^*)=0$ and
$ \hat{d}_{k+1}(x^*)=0,$ yielding
\begin{equation}
  \EXP{\|\hat{d}_{k+1}(x^*)\|^2 \mid F_{k}^{m}} = \hat{d}_{k}(x^*).
  \label{eqn:stu2}
\end{equation}
When $\hat{x}_k \notin L_{N},$ $\hat{x}_k = x_k$ and $\hat{x}_{k+1} =
x_{k+1}.$ Using relation (\ref{eqn:rel1}), we conclude that 
\begin{align*}
  \EXP{\|\hat{d}_{k+1}(x^*)\|^2 \mid F_{k}^{m} } 
  \leq& \|\hat{d}_k(x^*) \|^2 
   - 2 \alpha \left( f(\hat{x}_k) - f(x^*) \right) + 2m
  \alpha\mu \max_{x,y\in X} \|x-y\| \\
   & + \alpha^2 \left(\sum_{i=1}^{m} C_i+ m \nu \right)^2.
\end{align*}
Observe that when
$\hat{x}_k \notin L_N,$
\begin{align*}
f(\hat{x}_k) - f^* \geq \frac{1}{N} + m\mu \max_{x,y\in X} \|x-y\|
 + \frac{\alpha}{2}
\left(\sum_{i=1}^{m} C_i + m\nu \right)^2.
\end{align*}
Therefore, by combining the preceding two relations, we obtain 
for $\hat{x}_k \notin L_N,$
\begin{equation}
  \EXP{\|\hat{d}_{k+1}(x^*)\|^2 \mid F_{k}^{m} } 
  \le \|\hat{d}_k (x^*) \|^2 -\frac{2\alpha}{N}.
 \label{eqn:stu1}
\end{equation}
Therefore, from (\ref{eqn:stu2})~and~(\ref{eqn:stu1}), we can write
\begin{align}
  \EXP{\|\hat{d}_{k+1}(x^*)\|^2 \mid F_{k}^{m}} \leq \|\hat{d}_k (x^*)
  \|^2 - \Delta_{k+1},  
   \label{eqn:apply}
\end{align}
where 
\begin{align*}
  \Delta_{k+1} = \begin{cases} 0 &\mbox{ if $\hat{x}_k \in L_N$,} \\
  \frac{2 \alpha }{N} &\mbox{ if $\hat{x}_k \notin L_N.$}
  \end{cases}
\end{align*}
Observe that (\ref{eqn:apply}) satisfies the conditions of
Lemma~\ref{lemma:asm} with $u_{k} = \|\hat{d}_k (x^*) \|^2,$
$\mathcal{F}_k = F_{k}^{m},$ $q_{k} = 0,$ $v_{k} = \Delta_{k+1}$ and 
$w_{k}=0.$ Therefore, it follows that with
probability~1,
\[
\sum_{k=0}^{\infty} \Delta_{k+1} < \infty.
\]
However, this is possible only if $\Delta_{k} = 0$ for all $k$
sufficiently large. Therefore, with probability 1, we have $x_k\in
L_N$ for all sufficiently large $k$.  By letting $N \to \infty,$ we
obtain (\ref{eqn:result2b}).
\end{proof}

As seen from relation (\ref{eqn:result2b}) 
of Theorem \ref{thm:igac1},
the error bound on the ``best function'' value $\inf_k f(x_k)$ depends
on the step-size $\a$, and the bounds $\mu$ and $\nu$ for the moments
of the subgradient errors $\e_{i,k}$.  When the errors
$\e_{i,k}$ have zero mean, the results of Theorem \ref{thm:igac1} hold
with $\mu=0$. The resulting error bound is 
$\frac{\alpha}{2} \left(\sum_{i=1}^{m} C_i + m\nu\right)^2$,
which can be controlled with the step-size $\a$.
However, this result also holds when
the boundedness of $X$ is relaxed by requiring subgradient boundedness
instead, as seen in the following theorem.  The proof of this theorem
is similar to that of Theorem \ref{thm:igac1}, with some extra details
to account for the possibility that the optimal set $X^*$ may be empty.

\begin{theorem}
  \label{thm:igac2}
  Let Assumptions~\ref{ass:convex}, \ref{ass:erindep}, and
  \ref{ass:sgdbound} hold.  Let the sequence $\{x_k\}$ be generated by
  the method (\ref{eqn:stocincgrad}) with a constant step-size rule,
  i.e., $\alpha_k = \alpha$ for all $k \in \mathbb{N}$.  Also, assume
  that the subgradient errors $\e_{i,k}$ have zero mean and bounded
  second moments, i.e.,
  \[\mu_k=0\hbox{\ \ \ for all $k\ge1$},\qquad\quad
  \nu=\sup_{k \ge1}\nu_{k}<\infty.\]
  We then have
  \begin{equation}
   \liminf_{k\to\infty} \EXP{f(x_k)} \leq f^* 
     + \frac{\alpha}{2} \left(\sum_{i=1}^{m} C_i + m\nu\right)^2, 
   \label{eqn:result1b2}
   \end{equation}
and with probability 1,
   \begin{equation}
     \inf_{k\geq 0} f(x_k) \leq f^* 
      + \frac{\alpha}{2}
     \left(\sum_{i=1}^{m} C_i + m \nu \right)^2.
     \label{eqn:result2b2} 
  \end{equation}
\end{theorem}
\begin{proof}
All
the assumptions of Lemma \ref{lemma:key} are satisfied.
Since $\mu_k=0$  and $\nu_k\le \nu$ for 
all $k$,
according to the relation of
Lemma~\ref{lemma:key}, 
we have for any $y\in X$,
\begin{align}
    \EXP{\|d_{k+1}(y)\|^2 \mid F_{k}^{m} } \leq& \|d_{k}(y)\|^2 
    - 2\a\left(f(x_k) - f(y) \right)
    + \a^2\left(
    \sum_{i=1}^{m} C_i+ m \nu \right)^2.
  \label{eqn:rel2}
  \end{align}
By taking the total expectation, we obtain for all $y\in X$ and all $k$,
\begin{align}
  \EXP{\|d_{k+1}(y)\|^2} \leq& \EXP{\|d_k (y) \|^2} - 2 \alpha
  \left(\EXP{f(x_k)} - f(y) \right) 
  + \alpha^2 \left( \sum_{i=1}^{m} C_i+ m \nu \right).
  \label{eqn:exp02}
\end{align}
Assume now that the relation (\ref{eqn:result1b2}) is not valid.  Then
there will exist a $\gamma > 0$ and an index $k_{\gamma}$ such that
for all $k > k_{\gamma}$,
\[
\EXP{f(x_k)} \geq f^* + 2\gamma  + 
\frac{\alpha}{2} \left(\sum_{i=1}^{m} C_i + m\nu \right)^2.
\]
Let $y_\gamma\in X$ be such that $f(y_\gamma)\le f^*+\gamma$.
Therefore, for $k \ge k_{\gamma}$, we have
\[
\EXP{f(x_k)} -f(y_\gamma)\geq \gamma + \frac{\alpha}{2}
\left(\sum_{i=1}^{m} C_i + m\nu \right)^2.
\]
Fix $y = y_{\gamma}$ in (\ref{eqn:exp02}) and in a manner identical to
the proof of Theorem \ref{thm:igac1} we can obtain a contradiction.

To prove the relation in (\ref{eqn:result2b2}), we use a line of
analysis similar to that of the proof of Theorem \ref{thm:igac1},
where we define the set
\begin{align*}
L_{N} &= \left\{x \in X: f(x) < f^* + \frac{2}{N} + \frac{\alpha}{2}
       \left(\sum_{i=1}^{m} C_i + m\nu \right)^2 \right\}.
\end{align*}
We let  $y_N \in L_N$ be such that
\[f(y_N)\le f^*+\frac{1}{N},\]
and consider the sequence $\hat{x}_{k}$ defined as follows:
\[
\hat{x}_{k+1} = \begin{cases} x_{k+1} &\mbox{if  $\hat{x}_{k}
  \notin L_N$,} \\ y_N &\mbox{if $\hat{x}_{k} \in L_N$}. \end{cases}
\]
As in the proof of Theorem~\ref{thm:igac1}, we can show that the
sequence $\{x_{k}\}$ enters the set $L_N,$ for any $N.$ We then obtain
the result by taking the limit $N \to \infty.$
\end{proof}

In the absence of errors ($\nu=0$), the error bound of 
Theorem \ref{thm:igac2} reduces to 
\[f^* +\frac{\a}{2}\l(\sum_{i=1}^m C_i\r)^2,\] which coincides with the error 
bound for the cyclic incremental subgradient method (without errors)
established in \cite{Nedic01}, Proposition 2.1.

\section{Markov Randomized Incremental Subgradient Method}
\label{sec:markov}
While the method of Section \ref{sec:incgrad} is implementable in
networks with a ring structure (the agents form a cycle), the method
of this section is implementable in networks with an arbitrary
connectivity structure. The idea is to implement the incremental
algorithm by allowing agents to communicate only with their
neighbors. In particular, suppose at time $k$, an agent $j$ updates
and generates the estimate $x_k$. Then, agent $j$ may pass this
estimate to his neighboring agent $i$ with probability $[P(k)]_{i,j}$.
If agent $j$ is not a neighbor of $i$, then this probability is $0$.
Formally, the update rule for this method is given by
\begin{equation}
  x_{k+1} = \mathcal{P}_{X} \left[ x_{k} - \alpha_{k+1} \left(\nabla
  f_{s(k+1)}(x_{k}) +\e_{s(k+1),k+1} \right) \right], \label{eqn:mark}
\end{equation}
where $x_0\in X$ is some random initial vector,
$\epsilon_{s(k+1),k+1}$ is a random noise vector and $\alpha_{k+1}>0$
is the step-size. The sequence of indices of agents updating in time
evolves according to a a time non-homogeneous Markov chain with states
$1,\dots, m.$ We let $P(k)$ denote the transition matrix of this chain
at time $k,$ i.e.,
$$[P(k)]_{i,j}=\Prob\{s(k+1)=j\mid s(k)=i\}\qquad \hbox{for all }
i,j\in\{1,\ldots,m\}.$$ In the absence of subgradient errors
($\e_{s(k),k}=0$), when the probabilities $[P(k)]_{i,j}$ are all
equal, i.e., $[P(k)]_{i,j}=\frac{1}{m}$ for all $i,j$ and all $k$, the
method reduces to the incremental method with randomization proposed
in \cite{Nedic01}, which is applicable only to the agent networks that
are fully connected.

We note here that the time non-homogeneous Markov chain models
networks where the set of neighbors of an agent may change in time (as
the network may be mobile or for other reasons). We will also assume
that the agents decide the probabilities with which they communicate
with their neighbors, i.e., at time $k$, the agent $j$ chooses the
probabilities $[P(k)]_{i,j}\ge 0$ for his neighbors $i.$
 
The main difficulty in the analysis of the method in (\ref{eqn:mark})
comes from the dependence between the random agent index $s(k+1)$ and
the iterate $x_k$.  Assuming that the Markov chain is ergodic with the
uniform steady-state distribution, in the absence of the errors
$\e_{i,k}$ (i.e., $\e_{i,k}=0$), it is intuitively possible that the
method uses directions that approximate the subgradient $\frac{1}{m}
\sum_{i=1}^m \nabla f_{i}(x_{k})$ in the steady state.  This is the
basic underlying idea that we exploit in our analysis.

For this idea to work, it is crucial not only that the Markov chain
probabilities converge to a uniform distribution but also that the
convergence rate estimate is available in an explicit form.  The
uniform steady state requirement is natural since it corresponds to
each agent updating his objective $f_i$ with the same steady state
frequency, thus ensuring that the agents cooperatively minimize the
overall network objective function $f(x)=\sum_{i=1}^m f_i(x),$ and not
a weighted sum.  We use the rate estimate of the convergence of the
products $P(\ell)\cdots P(k)$ to determine the step-size choices that
guarantee the convergence of the method in (\ref{eqn:mark}).

To ensure the desired limiting behavior of the Markov chain
probabilities, we use the following two assumptions 
on the matrices $[P(k)]$.  
 
\begin{assumption}
  \label{ass:connect}
  Let $V = \{1,\ldots,m\}$. Let $E(k)$ be the set of edges $(j,i)$
  induced by the positive entries of the probability matrix $P(k)$,
  i.e.,
  $$E(k)=\{(i,j)\mid [P(k)]_{i,j} >0\}.$$ There exists an integer
  $Q\ge1$ such that the graph $\l(V,\cup_{l=k}^{k+Q-1} E(\ell)\r)$ is
  strongly connected for all $k.$
\end{assumption}

Generally speaking, Assumption \ref{ass:connect} ensures that each
agent has a ``chance'' to update the estimate $x_k$ once within a
finite time interval. It would guarantee that each agent updates the
estimate $x_k$ with the same frequency in a long run.  This is ensured
by the following assumption.

\begin{assumption}
  \label{ass:probab}
 \begin{itemize}
\item [(a)]
The diagonal entries of $P(k)$ are all positive for each $k$.
  \item[(b)] All positive entries of $[P(k)]$ are uniformly bounded
away from zero, i.e., there exists a scalar $\eta>0$ such that for all
$i,j\in\{1,\ldots,m\}$ and all $k$,
\[\hbox{if} \quad [P(k)]_{i,j}>0,
   \qquad\hbox{then}\quad [P(k)]_{i,j}>\eta.\]
  \item [(c)]
The matrix $P(k)$ is doubly stochastic for each $k$, 
i.e., the sum of the entries in every
    row and every column is equal to $1.$ 
  \end{itemize}
\end{assumption}

Assumptions \ref{ass:probab}(a) and (b) ensure that the information
from each and every agent is persistent in time.  Assumption
\ref{ass:probab}(c), ensures that the limiting Markov chain
probability distribution (if one exists) is uniform.  Assumptions
\ref{ass:connect} and \ref{ass:probab} together guarantee the
existence of the uniform limiting distribution, as shown in
\cite{Nedic08c}.  We state this result in the next section.

Note that the cyclic incremental algorithm (\ref{eqn:incgrad}) does
not satisfy Assumption~\ref{ass:probab}.  The transition probability
matrix corresponding to the cyclic incremental method is a permutation
matrix with the $(i,i)$-th entry being zero when agent $i$ updates at
time $k$. Thus, Assumption~\ref{ass:probab}(c) is violated.

We now provide some examples of transition matrices $[P(k)]$
satisfying Assumption \ref{ass:probab}. The second and third examples
are variations of the Metropolis-Hasting weights
\cite{Jadbabaie03,Xiao07}, defined in terms of the agent neighbors.
We let $N_i(k)\subset\{1,\ldots, m\}$ be the set of neighbors of an
agent $i$ at time $k$, and let $|N_i(k)|$ be the cardinality of this
set.  Consider the following rules:

\begin{itemize}
\item {\it Equal probability scheme}. The probabilities that agent
$i$ uses at time $k$ are 
\[[P(k)]_{i,j}=\begin{cases}
\frac{1}{m}
& \hbox{if } j\ne i \hbox{ and }j\in N_i(k), \cr \hbox{} &\hbox{}\cr
1-\frac{| N_i(k) |}{m}
& \hbox{if } j=i,\cr
0 & \hbox{otherwise}. \cr
 \end{cases}
\]
\item
{\it Min-equal neighbor scheme}. 
The probabilities that agent $i$ uses at time $k$ are 
\[[P(k)]_{i,j}=\begin{cases}
\min\l\{\frac{1}{|N_i(k)| +1},\frac{1}{|N_j(k)|+1}\r\} 
& \hbox{if } j\ne i \hbox{ and }j\in N_i(k),
\cr \hbox{} &\hbox{}\cr
1-\sum_{j\in N_i(k)} \min\l\{\frac{1}{|N_i(k)| +1},\frac{1}{|N_j(k)|+1}\r\}  
& \hbox{if } j=i,\cr
0 & \hbox{otherwise}. \cr
 \end{cases}
\]

\item
{\it Weighted Metropolis-Hastings scheme}.  The probabilities that
agent $i$ uses at time $k$ are given by
\[[P(k)]_{i,j}=\begin{cases}
\eta_i\,  \min\l\{\frac{1}{|N_i(k)|},\frac{1}{N_j(k)|}\r\} 
& \hbox{if } j\ne i \hbox{ and }j\in N_i(k),\cr
\hbox{} &\hbox{}\cr
1-\eta_i\,\sum_{j\in N_i(k)} \min\l\{\frac{1}{|N_i(k)|},\frac{1}{|N_j(k)|}\r\} 
& \hbox{if } j=i,\cr
0 & \hbox{otherwise}, \cr
 \end{cases}
\]
where the scalar $\eta_i>0$ is known only to agent $i$.
\end{itemize}

In the first example, the parameter $\eta$ can be defined as
$\eta=\frac{1}{m}.$ In the second example, $\eta$ can be defined as
$$\eta=\min_{i,j}\l\{\frac{1}{|N_i(k)| +1},\frac{1}{|N_j(k)|+1}\r\},$$
while in the third example, it can be defined as
$$\eta=\min_{i}\{\eta_i, 1-\eta_i\}\,
\min_{i,j}\l\{\frac{1}{|N_i(k)|},\frac{1}{|N_j(k)|}\r\}.$$ 

Furthermore, note that in the first example, each agent knows the size
of the network and no additional coordination with the other agents is
needed.  In the other two examples, an agent must be aware of the
number of the neighbors each of his neighbors has at any time.

\subsection{Preliminaries}
We first state a result from \cite{Nedic08d} for future reference. The
result captures the convergence and the rate of convergence of the
time non-homogeneous Markov chain to its steady state. Define
$\Phi(k,\ell),$ with $k>l,$ to be the transition probability matrix
for the Markov chain from time $l$ to $k$, i.e.\ $\Phi(k,\ell) =
P(\ell)\cdots P(k)$ with $k\ge l.$ Then, we have the following
convergence result for the transition matrices.
\begin{lemma}
  \label{lemma:roc}
  Assume the matrices $P(k)$ satisfy
  Assumptions~\ref{ass:connect} and \ref{ass:probab}. Then:
  \begin{enumerate}
  \item $\lim_{k \to \infty} \Phi(k,s) = \frac{1}{m} e e^T$ for all
    $s.$
  \item The convergence is geometric and the rate of
    convergence is given by
    \begin{align*}
      \left| [\Phi(k,\ell)]_{i,j} - \frac{1}{m} \right| \leq b\,
      \beta^{k-l}\qquad\hbox{for all $k$ and $l$ with $k\ge l\ge 0$},
    \end{align*}
    where 
    \begin{align*}
      b = \left(1 - \frac{\eta}{4m^2}\right)^{-2} \ \ \mbox{and} \ \ 
      \beta = \left(1 - \frac{\eta}{4m^2}\right)^{\frac{1}{Q}}.
    \end{align*}
  \end{enumerate}
\end{lemma}

We use the estimate of Lemma \ref{lemma:roc} to establish a key
relation in Lemma \ref{lemma:k2}, which is repeatedly invoked in our
subsequent analysis.  The idea behind Lemma \ref{lemma:k2} is the
observation that when there are no errors ($\e_{s(k),k}=0$) and the
Markov chain has a uniform steady state distribution, the directions
$\nabla f_{s(k+1)}(x_k)$ used in (\ref{eqn:mark}) are approximate
subgradients of the function $\frac{1}{m} \sum_{i=1}^m \nabla
f_{i}(x)$ at points $x_{n(k)}$ far away from $x_k$ in the past [i.e.,
$k>>n(k)$]. However, even though $x_n(k)$ are far away from $x_k$ in
time, their Euclidean distance $\|x_k-x_{n(k)}\|$ can be small when
the step-size is selected appropriately.  Overall, this means that
each iterate of method in (\ref{eqn:mark}) can be viewed as an
approximation of the iteration
$$x_{k+1}=P_X\left[x_k-\frac{\a_{k+1}}{m}\,\sum_{i=1}^m \nabla
f_i(x_k) + \a_{k+1}\xi_k\right],$$ with correlated errors $\xi_k$
depending on current and past iterates.

In the forthcoming lemma and thereafter, we let $G_{k}$ denote the
entire history of the method up to time $k$, i.e., the $\s$-field
generated by the initial vector $x_0$ and $\l\{ s(n),\e_{s(n),n}; 0
\leq n \leq k \r\}.$

\begin{lemma}
\label{lemma:k2}
  Let Assumptions~\ref{ass:convex}--\ref{ass:probab} hold.  Then, the
  iterates generated by algorithm (\ref{eqn:mark}) are such that for
  any step-size rule, for any $y\in X,$ and any non-negative integer
  sequence $\{n(k)\},$ $n(k) \leq k,$ we have
  \begin{align*}
  \EXP{\|d_{k+1}(y)\|^2 \mid G_{n(k)}} \leq& \EXP{\|d_{k}(y)\|^2 \mid
    G_{n(k)}} - \frac{2 \alpha_{k+1}}{m} \left(f\l(x_{n(k)}\r) - f(y)
    \right) \nonumber \\& + 2 b \l(\sum_{i=1}^{m} C_i\r) \alpha_{k+1}
    \beta^{k+1 - n(k)} \l\|d_{n(k)}(y)\r\| \nonumber \\& + 2 C
    \alpha_{k+1} \sum_{l = n(k)}^{k-1} \alpha_{l+1} \l( C +
    \nu_{l+1}\r) \nonumber \\&+ 2 \alpha_{k+1} \mu_{k+1} \EXP{\l\|
    d_{k}(y)\r\| \mid G_{n(k)}} + \alpha_{k+1}^2 (\nu_k + C)^2,
  \end{align*}
  where $d_k(y) = x_k - y$ and $C = \max_{1\le i\le m} C_i.$
\end{lemma}
\begin{proof}
Using the iterate update rule in (\ref{eqn:mark}) and the
non-expansive property of the Euclidean projection, we obtain for any
$y\in X$ and $k\ge0$,
\begin{align*}
    \|d_{k+1}(y)\|^2 =& \left\|\mathcal{P}_{X} \left[ x_{k} -
      \alpha_{k+1} \nabla f_{s(k+1)}\l(x_{k}\r) - \alpha_{k+1}
      \epsilon_{s(k+1),k+1} \right] - y \right\|^2 \nonumber \\ \leq&
      \left\| x_{k} - \alpha_{k+1} \nabla f_{s(k+1)}\l(x_{k}\r) -
      \alpha_{k+1} \epsilon_{s(k+1),k+1} - y \right\|^2 \nonumber \\
      =& \|d_{k}(y)\|^2 - 2 \alpha_{k+1} d_{k}(y)^T \nabla
      f_{s(k+1)}\l(x_{k}\r) \\& - 2 \alpha_{k+1} d_{k}(y)^T
      \epsilon_{s(k+1),k+1} + \alpha_{k+1}^2\left\|
      \epsilon_{s(k+1),k+1} + \nabla f_{s(k+1)}\l(x_{k}\r) \right\|^2.
\end{align*}
Using the subgradient inequality in (\ref{eqn:lip3}) to bound
$d_{k}(y)^T \nabla f_{s(k+1)}\l(x_{k}\r),$ we get
\begin{align*}
    \|d_{k+1}(y)\|^2 \leq& \|d_{k}(y)\|^2 - 2 \alpha_{k+1}
      \left(f_{s(k+1)}\l(x_{k}\r) - f_{s(k+1)}(y) \right) \\& - 2
      \alpha_{k+1} d_{k}(y)^T \epsilon_{s(k+1),k+1} +
      \alpha_{k+1}^2\left\| \epsilon_{s(k+1),k+1} + \nabla
      f_{s(k+1)}\l(x_{k}\r) \right\|^2 \\ =& \|d_{k}(y)\|^2 - 2
      \alpha_{k+1} \left(f_{s(k+1)}\l(x_{k}\r) -
      f_{s(k+1)}\l(x_{n(k)}\r) \right) \\& - 2 \alpha_{k+1}
      \left(f_{s(k+1)}\l(x_{n(k)}\r) - f_{s(k+1)}(y) \right) - 2
      \alpha_{k+1} d_{k}(y)^T \epsilon_{s(k+1),k+1} \\&+
      \alpha_{k+1}^2\left\| \epsilon_{s(k+1),k+1} + \nabla
      f_{s(k+1)}\l(x_{k}\r) \right\|^2.
\end{align*}
Taking conditional expectations with respect to the $\s$-field
$G_{n(k)}$, we obtain
\begin{align}
  \EXP{\|d_{k+1}(y)\|^2 \mid G_{n(k)}} \leq& \EXP{\|d_{k}(y)\|^2 \mid
      G_{n(k)}} \nonumber \\& - 2 \alpha_{k+1} \left(
      \EXP{f_{s(k+1)}\l(x_{k}\r) - f_{s(k+1)}\l(x_{n(k)}\r) \mid
      G_{n(k)}} \right) \nonumber\\&- 2 \alpha_{k+1} \left(
      \EXP{f_{s(k+1)}\l(x_{n(k)}\r) - f_{s(k+1)}\l(y\r) \mid G_{n(k)}}
      \right) \nonumber \\& - 2 \alpha_{k+1} \EXP{d_{k}(y)^T
      \epsilon_{s(k+1),k+1} \mid G_{n(k)}}\nonumber \\&+
      \alpha_{k+1}^2 \EXP{\left\| \epsilon_{s(k+1),k+1} + \nabla
      f_{s(k+1)}\l(x_{k}\r) \right\|^2 \mid
      G_{n(k)}}. \label{eqn:bound2}
\end{align}  
We next use the subgradient inequality in (\ref{eqn:lip3}) to estimate
the second term in the preceding relation.  
{\small
  \begin{align}
  \EXP{ f_{s(k+1)}\l(x_{k}\r) - f_{s(k+1)}\l(x_{n(k)}\r) \mid G_{n(k)}
  } \geq& \EXP{ \nabla f_{s(k+1)}\l(x_{n(k)}\r)^T \left( x_{n(k)} -
  x_{k} \right) \mid G_{n(k)} } \nonumber\\ \geq& - \EXP{ \l\|\nabla
  f_{s(k+1)}\l(x_{n(k)} \r)\r\| \left\| x_{n(k)} - x_{k} \right\| \mid
  G_{n(k)}} \nonumber\\\geq& -C \EXP{\left\| x_{n(k)} - x_{k} \right\|
  \mid G_{n(k)}}. \label{eqn:bound1}
\end{align}}
In the last step we have used the subgradient boundedness from
Assumption~\ref{ass:sgdbound} to bound the subgradient norms
$\l\|\nabla f_{s(k+1)}\l(x_{n(k)} \r)\r\|$ by $C=\max_{1\le i\le m}
C_i.$ We estimate $\EXP{\left\| x_{n(k)} - x_{k} \right\| \mid
G_{n(k)}}$ from the iterate update rule (\ref{eqn:mark}) and the
non-expansive property of the Euclidean projection as follows:
\begin{align*}
&\EXP{\l\| x_{n(k)} - x_{k} \r\| \mid G_{n(k)}} \\
&\ \ \ \leq \sum_{l= n(k)}^{k-1} 
\EXP{\l\| x_{l+1} - x_{l}\r\| \mid G_{n(k)} } \\ 
&\ \ \ \leq \sum_{l = n(k)}^{k-1}
\alpha_{l+1} \EXP{\l\|\nabla f_{s(\ell+1)}\l(x_{l}\r)\r\| +
\l\|\epsilon_{s(\ell+1),l+1}\r\| \mid G_{n(k)}} \\
&\ \ \ \leq \sum_{l =
n(k)}^{k-1} \alpha_{l+1} \EXP{\l\|\nabla f_{s(\ell+1)}\l(x_{l}\r)\r\| +
\EXP{\l\|\epsilon_{s(\ell+1),l+1}\r\| \mid G_{l} } \mid G_{n(k)}} \\
& \ \
\ \leq \sum_{l=n(k)}^{k-1} \alpha_{l+1} \left(C + \nu_{l+1} \right),
\end{align*}
where in the last step we have used the boundedness of subgradients
and of the second moments of $\e_{i,k}$ [cf.\ Eq.\
(\ref{eqn:jensen})].  From the preceding relation and Eq.\
(\ref{eqn:bound1}), we obtain
\begin{align*}
\EXP{ f_{s(k+1)}\l(x_{k}\r) - f_{s(k+1)}\l(x_{n(k)}\r) \mid G_{n(k)}}
  \geq& -C \sum_{l = n(k)}^{k-1} \alpha_{l+1} \l( C + \nu_{l+1}\r).
\end{align*}

By substituting the preceding estimate in (\ref{eqn:bound2}), we further 
obtain
\begin{align}
  \EXP{\|d_{k+1}(y)\|^2 \mid G_{n(k)}} \leq& \EXP{\|d_{k}(y)\|^2 \mid
      G_{n(k)}} + 2 C \alpha_{k+1} \sum_{l = n(k)}^{k-1} \alpha_{l+1}
      \l( C + \nu_{l+1}\r) \nonumber \\&- 2 \alpha_{k+1} \left(
      \EXP{f_{s(k+1)}\l(x_{n(k)}\r) - f_{s(k+1)}\l(y\r) \mid G_{n(k)}}
      \right) \nonumber \\& - 2 \alpha_{k+1} \EXP{d_{k}(y)^T
      \epsilon_{s(k+1),k+1} \mid G_{n(k)}}\nonumber \\&+
      \alpha_{k+1}^2 \EXP{\left\| \epsilon_{s(k+1),k+1} + \nabla
      f_{s(k+1)}\l(x_{k}\r) \right\|^2 \mid
      G_{n(k)}}. \label{eqn:status2}
\end{align}
We estimate the last term in (\ref{eqn:status2}) by using the
subgradient boundedness of Assumption~\ref{ass:sgdbound} and the
boundedness of the second moments of $\e_{i,k}$ [cf.\ Eq.\
(\ref{eqn:jensen})], as follows:
\begin{align*}
 &\EXP{\left\| \epsilon_{s(k+1),k+1} + \nabla f_{s(k+1)}\l(x_{k}\r)
      \right\|^2 \mid G_{n(k)}} \\
    & \ \ \ 
      \leq \EXP{\l\|
      \epsilon_{s(k+1),k+1} \r\|^2 + \left\|
      \nabla f_{s(k+1)}\l(x_{k}\r) \right\|^2  
      + 2 \l\| \epsilon_{s(k+1),k+1}\r\| \l\| \nabla
      f_{s(k+1)}\l(x_{k}\r) \right\|\mid G_{n(k)}} \\ 
      &\ \ \ \ 
      \leq\nu_k^2 + C^2 + 2 \nu_k C \\
      &\ \ \ \ 
      = (\nu_k + C)^2.
\end{align*}
Substituting the preceding estimate in Eq.\ (\ref{eqn:status2}), we have
\begin{align}
  \EXP{\|d_{k+1}(y)\|^2 \mid G_{n(k)}} \leq& \EXP{\|d_{k}(y)\|^2 \mid
    G_{n(k)}} + 2 C \alpha_{k+1} \sum_{l = n(k)}^{k-1} \alpha_{l+1}
    \l( C + \nu_{l+1}\r) \nonumber \\&- 2 \alpha_{k+1} \left(
    \EXP{f_{s(k+1)}\l(x_{n(k)}\r) - f_{s(k+1)}\l(y\r) \mid G_{n(k)}}
    \right) \nonumber \\& - 2 \alpha_{k+1} \EXP{d_{k}(y)^T
    \epsilon_{s(k+1),k+1} \mid G_{n(k)}}\nonumber \\&+ \alpha_{k+1}^2
    (\nu_k + C)^2.
  \label{eqn:status3}
\end{align}
We next estimate the term $\EXP{d_{k}(y)^T \epsilon_{s(k+1),k+1} \mid
G_{n(k)}}.$ Since $G_{n(k)} \subset G_{k}$ and $d_{k}(y)$ is
$G_k$-measurable
\begin{align*}
  \EXP{d_{k}(y)^T \epsilon_{s(k+1),k+1} \mid G_{n(k)}} =&
    \EXP{\EXP{d_{k}(y)^T \epsilon_{s(k+1),k+1} \mid G_{k}} \mid
    G_{n(k)}} \\ =& \EXP{d_{k}(y)^T \EXP{ \epsilon_{s(k+1),k+1} \mid
    G_{k}} \mid G_{n(k)}} \\ \geq& - \EXP{ \l\| d_{k}(y)\r\|
    \l\|\EXP{\epsilon_{s(k+1),k+1} \mid G_{k}}\r\| \mid G_{n(k)}}
    \\\geq& -\mu_{k+1} \EXP{\l\| d_{k}(y)\r\| \mid G_{n(k)}},
\end{align*}
where the first equality follows from the law of iterated
conditioning. Using the preceding estimate in (\ref{eqn:status3}), we
obtain
 \begin{align}
      \EXP{\|d_{k+1}(y)\|^2 \mid G_{n(k)}} \leq& \EXP{\|d_{k}(y)\|^2
	\mid G_{n(k)}} + 2 C \alpha_{k+1} \sum_{l = n(k)}^{k-1}
	\alpha_{l+1} \l( C + \nu_{l+1}\r) \nonumber \\&- 2
	\alpha_{k+1} \left( \EXP{f_{s(k+1)}\l(x_{n(k)}\r) -
	f_{s(k+1)}\l(y\r) \mid G_{n(k)}} \right) \nonumber \\& + 2
	\alpha_{k+1} \mu_{k+1} \EXP{\l\| d_{k}(y)\r\| \mid G_{n(k)}} +
	\alpha_{k+1}^2 (\nu_k + C)^2.
      \label{eqn:status4}
    \end{align}

Finally, we consider the term $\EXP{f_{s(k+1)}\l(x_{n(k)}\r) -
f_{s(k+1)}\l(y\r) \mid G_{n(k)}}$, and we use the fact that the
probability transition matrix for the Markov chain $\{s(k)\}$ from
time $n(k)$ to time $k$ is $\Phi(k+1,n(k))=P(n(k))\cdots P(k)$. We
have
\begin{align}
   & \EXP{f_{s(k+1)}\l(x_{n(k)}\r) - f_{s(k+1)}\l(y\r) \mid G_{n(k)}}
  \nonumber \\& \ = \sum_{i=1}^{m} \l[ \Phi(k+1,n(k))
  \right]_{s(n(k)),i} \left( f_i\l(x_{n(k)}\r) - f_i(y) \right)
  \nonumber \\& \ \geq \sum_{i=1}^{m} \frac{1}{m} \l(
  f_i\l(x_{n(k)}\r) - f_i(y) \right) - \sum_{i=1}^{m} \l|\l[
  \Phi(k+1,n(k)) \right]_{s(n(k)),i} -\frac{1}{m}\r| \l|
  f_i\l(x_{n(k)}\r) - f_i(y) \right| \nonumber \\& \ \geq \frac{1}{m}
  \l( f\l(x_{n(k)}\r) - f(y) \right) - b \beta^{k+1 - n(k)}
  \sum_{i=1}^{m} \l| f_i\l(x_{n(k)}\r) - f_i(y)
  \right|, \label{eqn:bound3}
\end{align}
where at the last step we have used Lemma~\ref{lemma:roc}.  Using the
subgradient inequality (\ref{eqn:lip3}), we further have
\begin{align}
  \l| f_i\l(x_{n(k)}\r) - f_i(y) \right| \leq& C_i \l\| x_{n(k)} - y
  \r\| = C_i \| d_{n(k)}(y) \|. \label{eqn:bound4}
\end{align}
The result now follows by combining the relations in 
Eqs.\ (\ref{eqn:status4}), (\ref{eqn:bound3}) and (\ref{eqn:bound4}).
\end{proof}

\subsection{Convergence for diminishing step-size}
In this section, we establish the convergence of the Markov randomized
method in (\ref{eqn:mark}) for a diminishing step-size.  Recall that
in Theorem~\ref{thm:QISGA1} for the cyclic incremental method, we
showed an almost sure convergence result for a diminishing step-size
$\alpha_k$ subject to some conditions that coordinate the choice of
the step-size, and the bounds $\m_k$ and $\nu_k$ on the moments of the
errors $\e_{i,k}$.  To obtain an analogous result for the Markov
randomized method, we use boundedness of the set $X$ and more
restricted step-size. In particular, we consider a step-size of the
form $\alpha_k=\frac{a}{k^p}$ for a range of values of $p,$ as seen in
the following.
\begin{theorem}
  \label{thm:MISGA1}
  Let Assumptions~\ref{ass:convex}, \ref{ass:erindep},
  \ref{ass:connect} and \ref{ass:probab} hold.  Assume that the
  step-size is $\a_k = \frac{a}{k^p},$ where $a$ and $p$ are positive
  scalars with $\frac{2}{3} < p\le1.$ In addition, assume that the
  bounds $\m_k$ and $\nu_k$ on the error moments satisfy 
  \[\sum_{k=1}^{\infty} \alpha_k\mu_k < \infty,\qquad
    \nu = \sup_{k \ge1} \nu_{k} < \infty.   
  \]
  Furthermore, let the set $X$ be bounded.  Then, with probability 1,
  we have
  $$\liminf_{k\to\infty} f(x_k) = f^*,\qquad \liminf_{k\to\infty}{\rm
    dist}(x_k,X^*) = 0.$$ 
\end{theorem}
\begin{proof}
Since the set $X$ is compact and $f_i$ is convex over $\rn$, it
follows that the subgradients of $f_i$ are bounded over $X$ for each
$i$.  Thus, Assumption \ref{ass:sgdbound} is satisfied, and we can use
Lemma~\ref{lemma:k2}.

Since $X$ is compact and $f$ is convex over $\rn$ (therefore, also
continuous), the optimal set $X^*$ is nonempty, closed and convex.
Let $x^*_k$ be the projection of $x_k$ on the set $X^*$.  In
Lemma~\ref{lemma:k2}, we let $y=x_k^*$ and let $n(k) = k +1 -\l\lceil
k^{\gamma} \r\rceil,$ where $\gamma>0$ (to be specified more precisely
later on). Note that $n(k) \le k$ for all $k\ge 1$.  Using this and
the relation $\di{x_{k+1}}\le \|x_{k+1}-x_k^*\|$, from
Lemma~\ref{lemma:k2}, we obtain for all $k>1$\footnote{The equivalent
expression for the case when $k=1$ is obtain by setting the fourth
term to 0.},
\begin{align*}
  \EXP{\di{x_{k+1}}^2 \mid G_{n(k)}} \leq& \EXP{\di{x_k}^2
    \mid G_{n(k)}} - \frac{2 \alpha_{k+1}}{m} \left(f\l(x_{n(k)}\r) -
    f^* \right) \nonumber \\
    & + 2 b\, \l(\sum_{i=1}^{m} C_i\r) \alpha_{k+1}
    \beta^{\l\lceil k^{\gamma} \r\rceil} \l\|d_{n(k)}(x_k^*)\r\| \nonumber
    \\
    & + 2 C \alpha_{k+1} \alpha_{n(k)+1} (\l\lceil k^{\gamma}
    \r\rceil - 2) \max_{n(k)\leq l\le k}\l( C + \nu_{l+1}\r) \nonumber
    \\
    &+ 2 \alpha_{k+1} \mu_{k+1} \EXP{\di{x_k}\mid
    G_{n(k)}} + \alpha_{k+1}^2 (\nu_k + C)^2.
\end{align*}
Taking expectations and using $\sup_{k\ge1}\nu_k=\nu$, 
we obtain for all $k > 1$,
\begin{align*}
  \EXP{\di{x_{k+1}}^2} \leq \EXP{\di{x_k}^2} - \frac{2
      \alpha_{k+1}}{m} \left( \EXP{f\l(x_{n(k)}\r)} - f^* \right) +
      \tau_{k+1},
\end{align*}
where 
\begin{align*}
  \tau_{k+1} =& 2 b\, \l(\sum_{i=1}^{m} C_i\r) \alpha_{k+1}
  \beta^{\l\lceil k^{\gamma} \r\rceil} \l\|d_{n(k)}(x_k^*)\r\| \\
  &+ 2 C (C+\nu)
  \alpha_{k+1} \alpha_{n(k)+1} (\l\lceil k^{\gamma} \r\rceil - 2)
  \nonumber \\
   &+ 2
  \alpha_{k+1} \mu_{k+1} \EXP{\di{x_k}} +
  \alpha_{k+1}^2 (\nu_k + C)^2.
\end{align*}

We next show that $\sum_{k=2}^{\infty} \tau_{k+1}<\infty$.  Since
$\alpha_k = \frac{a}{k^p}$, we have $\a_{k+1}<\a_k$ for all $k\ge1.$
Furthermore, since $\beta<1$, we have $\beta^{\l\lceil
k^{\gamma}\r\rceil}< \beta^{k^{\gamma}}.$ Therefore, $\alpha_{k+1}
\beta^{\l\lceil k^{\gamma}\r\rceil} < \frac{a \beta^{k^{\gamma}}
}{k^p}.$ By choosing $\gamma>0$ such that $\gamma\ge 1-p$, we see that
$\frac{1}{k^p}\le \frac{1}{k^{1-\gamma}}$ for all $k> 1.$ Hence, for
all $k> 1$,
\begin{align*}
\sum_{k=2}^{\infty} \alpha_{k+1} \beta^{\l\lceil k^{\gamma}\r\rceil} <
\sum_{k=2}^{\infty} \frac{a \beta^{k^{\gamma}} }{k^p} \leq
\sum_{k=2}^{\infty} \frac{a \beta^{k^{\gamma}}}{k^{1-\gamma}} \leq a
\int_{1}^{\infty} \frac{\beta^{y^{\gamma}}}{y^{1-\gamma}}\ dy =
- \frac{a \beta}{\gamma \ln(\beta)}.
\end{align*}
Since the set $X$ is bounded, it follows that
\begin{equation}
\sum_{k=2}^{\infty} 
2 b\, \l(\sum_{i=1}^{m} C_i\r) \alpha_{k+1}
  \beta^{\l\lceil k^{\gamma} \r\rceil} \l\|d_{n(k)}(x^*_k)\r\|
< \infty.
\label{eqn:relone}
\end{equation}

Next, since $\l\lceil k^{\gamma} \r\rceil - 2<k^\gamma$ for all $k \ge
2,$ and since $\a_{k+1}<\a_k$, $\a_k=\frac{1}{k^p}$ and
$n(k)=k+1-\l\lceil k^{\gamma}\r\rceil$, it follows that for all
$k\ge2,$
\[\alpha_{k+1} \alpha_{n(k)+1} (\l\lceil k^{\gamma} \r\rceil - 2) 
< \frac{a^2\,k^{\gamma}}{k^p ( k + 2 - \l\lceil k^{\gamma} \r\rceil)^p}
<\frac{a^2\,k^{\gamma}}{k^p ( k - k^{\gamma})^p}
=\frac{a^2\,k^{\gamma}}{k^{2p} (1 - k^{\gamma-1})^p}.\]
By choosing $\gamma>0$ such that it also satisfies $\gamma <2p-1$
(in addition to $\gamma \ge 1-p$), we have $\gamma<1$ (in view of $p\le 1$).
Therefore, for all $k\ge 2$, 
\[\frac{k^{\gamma}}{k^{2p} (1 - k^{\gamma-1})^p}
\le \frac{1}{(1 - 2^{\gamma - 1})^p}\frac{1}{k^{2p - \gamma}}.\]
By combining the preceding two relations, we have 
{\small 
 \begin{align}
\sum_{k=2}^{\infty} 
  2C(C+\nu)\alpha_{k+1} \alpha_{n(k)+1} (\l\lceil
  k^{\gamma} \r\rceil - 2) 
  <  
   2C(C+\nu)\frac{a^2}{(1 - 2^{\gamma - 1})^p}
   \sum_{k=2}^{\infty} \frac{1}{k^{2p - \gamma}} <\infty,
  \label{eqn:reltwo}
\end{align}
} where the finiteness of the last sum follows from $2p-\gamma>1$.  

Finally, as a consequence of our assumptions, we also have
\begin{align*}
  &\sum_{k=2}^{\infty} 2 \alpha_{k+1} \mu_{k+1} \EXP{\di{x_k}} < \infty, \\ 
  &\sum_{k=2}^{\infty} \alpha_{k+1}^2
    (\nu_k + C)^2 < \infty.
\end{align*}
Thus, from Eqs.\ (\ref{eqn:relone}) and (\ref{eqn:reltwo}), and the
preceding two relations, we see that $\sum_{k=2}^{\infty} \tau_{k+1} <
\infty.$

From the deterministic analog of Lemma~\ref{lemma:asm} we conclude that
$\EXP{\di{x_k}^2}$ converges to a non-negative scalar and
\begin{align*}
  \sum_{k=2}^{\infty} \frac{2 \alpha_{k+1}}{m} \left(
  \EXP{f\l(x_{n(k)}\r)} - f^* \right) < \infty.
\end{align*}
Since $p<1,$ we have $\sum_{k=2}^{\infty} \a_{k+1}=\infty$.  Further,
since $f\l(x_{n(k)}\r) \geq f^*,$ it follows that
\begin{equation}
  \liminf_{k\to\infty} \EXP{f\l(x_{n(k)}\r)} = f^*.\label{eqn:liminfef}
\end{equation}
The function $f$ is convex over $\rn$ and, hence, continuous. Since
the set $X$ is bounded, the function $f(x)$ is also bounded on $X.$
Therefore, from Fatou's lemma it follows that
\[
\EXP{\liminf_{k\to\infty} f\l(x_{k}\r) } 
\leq \liminf_{k\to\infty} \EXP{f\l(x_{k}\r)}
= f^*,
\]
implying that $\liminf_{k\to\infty} f\l(x_{k}\r) = f^*$ with
probability 1.  Moreover, from this relation, by the continuity of $f$
and boundedness of $X$, it follows that $\liminf_{k\to\infty}
\di{x_{k}}= 0$ with probability 1.
\end{proof}

As seen in the proof of Theorem \ref{thm:MISGA1}, $\EXP{\di{x_{k}}^2}$
converges to a nonnegative scalar, but we have no guarantee that its
limit is zero.  However, this can be shown, for example, for a
function with a sharp set of minima, i.e., $f$ satisfying
\[
  f(x) - f^* \geq \zeta\, \di{x}^{\xi}\qquad\hbox{for all }x\in X,
  \]
  for some positive scalars $\zeta$ and $\xi.$ Under the assumptions
of Theorem \ref{thm:MISGA1}, we have that $\liminf_{k\to\infty}
\EXP{f\l(x_{k}\r)}= f^*$ [cf.\ Eq.\ (\ref{eqn:liminfef})] and
therefore,
\begin{align*}
0 = \liminf_{k\to\infty} \EXP{f\l(x_{k}\r) - f^*} 
\geq \zeta\, \liminf_{k\to\infty}
\EXP{\di{x_k}^{\xi}} \geq 0.
\end{align*}
Hence, $\liminf_{k\to\infty} \EXP{\di{x_{k}}^\xi} = 0,$ and since
 $\EXP{\di{x_k}^2}$ converges, it has to
converge to $0.$

\subsection{Error bounds for constant step-size}
We now establish error bounds when the Markov randomized incremental
method is used with a constant step-size.

\begin{theorem}
  \label{thm:dgac1}
  Let Assumptions~\ref{ass:convex}, \ref{ass:erindep},
  \ref{ass:connect} and \ref{ass:probab} hold. 
  Let the sequence $\{x_k\}$ be generated by
  the method (\ref{eqn:mark}) with a constant step-size rule, i.e.,
  $\alpha_k = \alpha$ for all $k$.  Also, assume that the set $X$ is
  bounded, and
  \[\mu=\sup_{k \ge1}\mu_{k}<\infty,\qquad 
    \nu=\sup_{k \ge1}\nu_{k}<\infty.\]
  Then for any integer $T \ge0,$
  \begin{align}
    \lim \inf_{k} \EXP{f(x_k)} \leq& f^* + \mu \max_{x,y\in X}\|x-y\| +
    \frac{1}{2} \alpha (\nu + C)^2 + \alpha T C \l( C + \nu
    \r)\nonumber\\
    & + b \l(\sum_{i=1}^{m} C_i\r) \beta^{T+1}
    \max_{x,y\in X}\|x-y\|, \label{eqn:lims}
  \end{align}
  where $\beta=\left(1-\frac{\eta}{4m^2}\right)^{\frac{1}{Q}}$ and $C
  = \max_{1\le i \le m} C_i.$ Furthermore, with probability 1, the
  same estimate holds for $\inf_{k} f(x_k)$.
  \vskip 1pc
\end{theorem}
\begin{proof}
Since $X$ is compact and each $f_i$ is convex over $\rn$, the
subgradients of $f_i$ are bounded over $X$ for each $i$. Thus, all the
assumptions of Lemma \ref{lemma:k2} are satisfied.  Let $T$ be a
nonnegative integer and let $n(k) = k -T$.  Since $\mu_k\le \mu$ and
$\nu_k\le \nu$ for all $k$, and $\| d_{k}(y) \|\le \max_{x,y\in
X}\|x-y\|,$ according to Lemma~\ref{lemma:k2}, we have
for $y=x^*\in X^*$ and $k \geq T,$
\begin{align}
  \EXP{\|d_{k+1}(x^*)\|^2 \mid G_{n(k)}} \leq& \EXP{\|d_{k}(x^*)\|^2
    \mid G_{n(k)}} - \frac{2 \alpha}{m} \left(f\l(x_{k - T}\r) - f^*
    \right) \nonumber \\
    & + 2 b\, \l(\sum_{i=1}^{m} C_i\r) \alpha
    \beta^{T+1} \max_{x,y\in X}\|x-y\| \nonumber \\
    & + 2 \alpha^2 T C
    \l( C + \nu \r) \nonumber \\
    &+ 2 \alpha \mu \max_{x,y\in X}\|x-y\|
    + \alpha^2 (\nu + C)^2.
    \label{eqn:dstep1}
\end{align}
By taking the total expectation, we obtain for all $x^*\in X^*$ and
all $k\ge T$,
\begin{align*}
   \EXP{\|d_{k+1}(x^*)\|^2} \leq& \EXP{\|d_{k}(x^*)\|^2} - \frac{2
    \alpha}{m} \left(\EXP{f\l(x_{k - T}\r)} - f^* \right) \nonumber\\
    & + 2 b\, \l(\sum_{i=1}^{m} C_i\r) \alpha \beta^{T+1} \max_{x,y\in
    X}\|x-y\| \nonumber \\
    & + 2 \alpha^2 T C \l( C + \nu \r) \nonumber\\
    &+ 2 \alpha \mu \max_{x,y\in X}\|x-y\| + \alpha^2 (\nu + C)^2.
\end{align*}
Now assume that the relation (\ref{eqn:lims}) does not hold.  Then,
there will exist a $\gamma > 0$ and an index $k_{\gamma}\ge T$ such
that for all $k \ge k_{\gamma}$,
\begin{align*}
\EXP{f(x_k)} \geq& f^* + \gamma + \mu \max_{x,y\in X}\|x-y\| +
    \frac{1}{2} \alpha (\nu + C)^2 + \alpha T C \l( C + \nu \r)\\
    & + b\,
    \l(\sum_{i=1}^{m} C_i\r) \beta^{T+1} \max_{x,y\in X}\|x-y\|.
\end{align*}
Therefore, for $k \ge k_{\gamma}+T$, we have
\begin{align*}
  \EXP{\|d_{k+1} (x^*) \|^2 } \leq& \EXP{\|d_{k} (x^*) \|^2 } - 2
  \alpha \gamma\le\cdots\le \EXP{\|d_{k_\gamma}(x^*)\|^2 } - 2
  \a\gamma (k - k_{\gamma}).
\end{align*}
For sufficiently large $k,$ the right hand side of the preceding
relation is negative, yielding a contradiction.  Thus, the relation
(\ref{eqn:lims}) must hold for all $T\ge0.$

We next show that for any $T\ge0,$
\begin{align}
\inf_{k} f(x_k) \leq & f^* + \mu \max_{x,y\in X}\|x-y\| + \frac{1}{2}
    \alpha (\nu + C)^2 + \alpha T C \l( C + \nu \r)\nonumber\\& + b\,
    \l(\sum_{i=1}^{m} C_i\r) \beta^{T+1} \max_{x,y\in X}\|x-y\|,
    \label{eqn:dresult2b}
\end{align}
with probability $1.$ Define the set
\begin{align*}
L_{N} = \left\{x \in X: f(x) < f^* + \frac{1}{N} + \mu \max_{x,y\in
    X}\|x-y\| + \frac{1}{2} \alpha (\nu + C)^2 + \alpha T C \l( C + \nu
    \r) \r.& \nonumber\\ \l. + b\, \l(\sum_{i=1}^{m} C_i\r) \beta^{T+1}
    \max_{x,y\in X}\|x-y\| \r\}.
\end{align*}
Let $x^* \in X^*$ and define the sequence $\hat{x}_{k}$ as follows:
\[
\hat{x}_{k+1} = \begin{cases} x_{k+1} &\mbox{if $\hat x_{k} \notin L_N$}, 
  \\ 
  x^* &\mbox{otherwise}. \end{cases}
\]
Thus, the process $\{\hat{x}_k\}$ is identical to the process
$\{x_k\}$ until  $\{x_k\}$ enters the set $L_N.$
Define 
\[
\hat{d}_k(y) = \hat{x}_k - y\qquad \hbox{ for any $y\in X$}.
\]
Let $k \geq T.$ Consider the case when $\hat{x}_{k} \in L_{N}.$ Then,
$\hat x_k=x^*$ and $\hat x_{k+1}=x^*$, so that $\hat{d}_{k}(x^*)=0$
and $\hat{d}_{k+1}(x^*)=0,$ yielding 
\begin{equation}
\EXP{\|\hat{d}_{k+1}(x^*)\|^2 \mid G_{n(k)}}
=\EXP{\|\hat{d}_{k}(x^*)\|^2 \mid G_{n(k)}}.
  \label{eqn:dstu2}
\end{equation}
Consider now the case when $\hat{x}_{k} \notin L_{N}.$ Then,
$\hat x_l=x_l$ and $x_l\notin L_N$ for all $l\le k+1$.
Therefore, by the definition of the set $L_N$, we have
\begin{align}
  f(x_{k - T})-f^* \geq 
   & \frac{1}{N} + \mu \max_{x,y\in X}\|x-y\|
    + \frac{1}{2} \alpha (\nu + C)^2 + \alpha T C \l( C + \nu \r) &
    \nonumber\\ 
   &+ b\, \l(\sum_{i=1}^{m} C_i\r) \beta^{T+1} \max_{x,y\in
    X}\|x-y\|. \label{eqn:dd}
\end{align}
By using relations (\ref{eqn:dstep1}) and (\ref{eqn:dd}), we conclude
that for $\hat{x}_k \notin L_{N},$
\begin{equation}
  \EXP{\|\hat{d}_{k+1}(x^*)\|^2 \mid G_{n(k)}} 
  \le \EXP{\|\hat{d}_k (x^*) \|^2 \mid G_{n(k)}} -\frac{2\alpha}{N}.
  \label{eqn:dstu1}
\end{equation}
Therefore, from (\ref{eqn:dstu2})~and~(\ref{eqn:dstu1}), we can write
\begin{align}
  \EXP{\|\hat{d}_{k+1}(x^*)\|^2 \mid G_{n(k)}} \leq
  \EXP{\|\hat{d}_k (x^*) \|^2 \mid G_{n(k)}} - 
  \Delta_{k+1},
   \label{eqn:dapply}
\end{align}
where 
\begin{align*}
  \Delta_{k+1} = \begin{cases} 0 &\mbox{ if $\hat{x}_k \in L_N,$} \\
  \frac{2 \alpha}{N} &\mbox{ if $\hat{x}_k \notin L_N.$}
  \end{cases}
\end{align*}
Observe that (\ref{eqn:dapply}) satisfies the conditions of
Lemma~\ref{lemma:asm} with $u_{k} = \EXP{\|\hat{d}_k (x^*) \|^2 \mid
G_{n(k)}},$ $\mathcal{F}_k = G_{n(k)},$ $q_{k} = 0,$ $w_{k} = 2 \Delta_{k+1}$ 
and $v_{k} = 0.$ Thus, it follows that with probability~1,
\[
\sum_{k=T}^{\infty} \Delta_{k+1} < \infty.
\]
However, this is possible only if $\Delta_{k} = 0$ for all $k$
sufficiently large. Therefore, with probability 1, we have $x_k\in
L_N$ for all sufficiently large $k$.  By letting $N \to \infty,$ we
obtain (\ref{eqn:dresult2b}).
\end{proof}

Under Assumptions of Theorem \ref{thm:dgac1}, the function $f$ is
bounded over the set $X$, and by Fatou's lemma, we have
\[\EXP{\liminf_{k\to\infty} f(x_k)}\le \liminf_{k\to\infty} \EXP{f(x_k)}.\]
It follows that the estimate of Theorem \ref{thm:dgac1} also holds for
$\EXP{\liminf_{k\to\infty} f(x_k)}.$

In the absence of errors ($\m_k=0$ and $\nu_k=0$),
the error bound in Theorem \ref{thm:dgac1}  reduces to 
\begin{equation}
f^* + \frac{1}{2}\, \a C^2 + \a TC^2 +b\l(\sum_{i=1}^m C_i\r)\beta^{T+1}
\max_{x,y\in X}\|x-y\|.\label{eqn:errbound}
\end{equation}
With respect to the parameter $\beta$, the error bound is obviously
smallest when $\beta=0$.  This corresponds to uniform transition
matrices $P(k)$, i.e., $P(k)=\frac{1}{m}ee^T$ for all $k$ (see Lemma
\ref{lemma:roc}).  As mentioned, the Markov randomized method with
uniform transition probability matrices $P(k)$ reduces to the
incremental method with randomization in \cite{Nedic01}.  In this
case, choosing $T=0$ in (\ref{eqn:errbound}) is optimal and the
resulting bound is $f^*+\frac{\a}{2}C^2$, with $C=\max_{1\le i\le m}
C_i$.  We note that this bound is better by a factor of $m$ than the
corresponding bound for the incremental method with randomization
given in Proposition 3.1 in\cite{Nedic01}.

When transition matrices are non-uniform ($\beta>0$), and good
estimates of the bounds $C_i$ on subgradient norms and the diameter of
the set $X$ are available, one may optimize the error bound in
(\ref{eqn:errbound}) with respect to integer $T$ for $T\ge0$. In
particular, one may optimize the term $ \a TC^2 +b\l(\sum_{i=1}^m
C_i\r)\beta^{T+1} \max_{x,y\in X}\|x-y\|$ over integers $T\ge0$.  It
can be seen that the optimal integer $T^*$ is given by
\begin{equation}
T^*= \begin{cases}
  0 &\hbox{when \ } \frac{\a C^2}{C_0 (-\ln
  \beta)}\ge 1,\cr \l\lceil \l(\ln \beta\r)^{-1}\, \ln \l(\frac{\a
  C^2}{C_0 (-\ln \beta)}\r) \r\rceil-1 & \hbox{when \ }\quad \frac{\a
  C^2}{C_0 (-\ln \beta)}< 1,
\end{cases}
\label{eqn:optimalT}
\end{equation}
where
$C_0 =b\,\l(\sum_{i=1}^m C_i\r) \max_{x,y\in X}\|x-y\|$.

A similar expression for optimal $T^*$ in the presence of subgradient
errors can be obtained, but it is rather cumbersome. Furthermore such
an expression (as well as the preceding one) may not be of practical
importance when the bounds $C_i$, the diameter of the set $X$, and the
bounds $\m$ and $\nu$ on the error moments are ``roughly'' known. In
this case, a simpler bound can be obtained by just comparing the
values $\a$ and $\beta$, as given in the following.

\begin{corollary}
  \label{cor:dgac1}
  Let the conditions of Theorem \ref{thm:dgac1} hold. Then,
  \begin{align*}
    \liminf_{k\to\infty} \EXP{f(x_k)} \leq& 
    	f^* + \mu \max_{x,y\in X}\|x-y\| \cr
            & + \alpha \left[\frac{1}{2} (\nu + C)^2 
                    + b\, \l(\sum_{i=1}^{m} C_i\r) 
                          \max_{x,y\in X}\|x-y\|\right]
            +\delta(\a,\beta),
  \end{align*} 
  where 
  \[\delta(\a,\beta)=
  \begin{cases}
    0 & \hbox{if } \alpha \geq \beta,\cr
    \l\lceil \frac{\ln(\alpha)}{\ln(\beta)}\r\rceil -1
    & \mbox{if }\alpha < \beta.
  \end{cases}
  \]  
  Furthermore, with probability 1, the same estimate holds for
  $\inf_{k} f(x_k)$.
\end{corollary}
\begin{proof}
When $\alpha > \beta$ choose $T = 0.$ In this case, from
(Theorem~\ref{thm:dgac1}) we get
\begin{align*}
\EXP{f(x_k)} \leq& f^* + \mu \max_{x,y\in X}\|x-y\| + \alpha \left(
    \frac{1}{2} (\nu + C)^2 + b\, \l(\sum_{i=1}^{m} C_i\r) \max_{x,y\in
    X}\|x-y\|\right).
\end{align*}
When $\alpha < \beta$ we can choose $T = \l\lceil
\frac{\ln(\alpha)}{\ln(\beta)}\r\rceil - 1.$ Then, from
(Theorem~\ref{thm:dgac1}), {\small
\begin{align*}
\EXP{f(x_k)} \geq 
    & f^* + \mu \max_{x,y\in X}\|x-y\|  \cr
    &+ \alpha \l[
    \frac{1}{2} (\nu + C)^2 + C \l( C + \nu \r) \l( \l\lceil
    \frac{\ln(\alpha)}{\ln(\beta)}\r\rceil - 1 \r) 
     + b\, \l(\sum_{i=1}^{m} C_i\r) \max_{x,y\in X}\|x-y\|\r].&
\end{align*}
}
\end{proof}

It can be seen that the error bounds in (\ref{eqn:optimalT}) and
Corollary~\ref{cor:dgac1} converge to zero as $\alpha\to0$. This is
not surprising in view of the convergence of the method with a
diminishing step-size.

As discussed earlier, the error bound in \cite{Johansson07} is
obtained assuming that there are no errors in subgradient evaluations
and that the sequence of computing agents form a homogeneous Markov
chain.  Here, while we relax these assumptions, we make the additional
assumption that the set $X$ is bounded.

A direct comparison between the bound in Corollary~(\ref{cor:dgac1})
and the results in \cite{Johansson07} is not possible.  However, some
qualitative comparisons on the nature of the bounds can be made. The
bound in \cite{Johansson07} is obtained for each individual agent's
sequence of iterates (by sampling the iterates).  This is a stronger
result than our results in (\ref{eqn:optimalT}) and
Corollary~\ref{cor:dgac1}, which provide guarantees only on the entire
iterate sequence (and not on the sequence of iterates at an individual
agent).  However, the bound in \cite{Johansson07} depends on the
entire network topology, through the probability transition matrix $P$
of the Markov chain. Thus, the bound can be evaluated \emph{only} when
the complete network topology is available.  In contrast, our bounds
given in (\ref{eqn:optimalT}) and Corollary~\ref{cor:dgac1} can be
evaluated without knowing the network topology.  We require that the
topology satisfies a connectivity assumption, as specified by
Assumption~\ref{ass:connect}, but we do not assume the knowledge of
the exact network topology.

\section{Discussion}
\label{sec:discuss}
Incremental algorithms form the middle ground between selfish agent
behavior and complete network cooperation. Each agent can be viewed to
be selfish, as it adjusts the iterate only using its own cost
function. At the same time, the agents also cooperate by passing the
iterate to a neighbor so that he may factor in his opinion by
adjusting the iterate using his cost function.  Through
Theorems~\ref{thm:QISGA1} and \ref{thm:MISGA1}, it was observed that a
system level global optimum could still be obtained through some
amount of cooperation. This can be construed as a statement of Adam
Smith's invisible hand hypothesis in a more semi-cooperative market
setting.

The results we have obtained are asymptotic in nature. The key step in
dealing with both the incremental algorithms was to obtain the basic
iterate equation (Lemmas~\ref{lemma:key} and \ref{lemma:k2}). This was
then combined with standard stochastic analysis techniques to obtain
asymptotic results. While we have restricted ourselves to establishing
only convergence results, it is possible to combine the techniques in
\cite{Nedic01b} with the basic iterate relation to obtain bounds on
the expected rate of convergence of the algorithms.  Finally, we have
only listed a few possible applications for the results in this
paper. The problem of aligning and coordinating mobile agents can also
be cast in the optimization framework studied in this paper and the
results in this paper, especially the results on Markov stochastic
sub-gradient algorithm, can be used to design suitable alignment
algorithms.

\bibliographystyle{siam} 
\bibliography{sism}
\end{document}